\documentclass[aps,prx,reprint,floatfix,superscriptaddress,amsmath,amssymb]{revtex4-2}
\usepackage[T1]{fontenc}
\usepackage[utf8]{inputenc}
\usepackage{mathtools}
\usepackage{bm}
\usepackage{dsfont}
\usepackage{subcaption}
\usepackage{caption}
\usepackage{bigints}
\usepackage{latexsym}
\usepackage[nointegrals]{wasysym}

\usepackage{graphicx}
\graphicspath{{figures/}}
\usepackage[dvipsnames]{xcolor}
\usepackage{multirow}
\usepackage{booktabs}
\usepackage{tcolorbox}
\usepackage{float}
\usepackage{algpseudocode}
\let\oldReturn\Return
\renewcommand{\Return}{\State\oldReturn}
\usepackage{enumerate}
\usepackage{beramono}
\usepackage{listings}
\definecolor{mustard}{HTML}{707238}
\definecolor{PaleGreen}{HTML}{98FB98}
\definecolor{LightRed}{HTML}{FE7A7A}
\definecolor{LightYellow}{HTML}{EFFA86}
\definecolor{YellowOrange}{HTML}{FFAE42}
\definecolor{OliveGreen}{HTML}{408000}
\definecolor{dkgreen}{rgb}{0,0.6,0}
\definecolor{gray}{rgb}{0.5,0.5,0.5}
\definecolor{mauve}{rgb}{0.58,0,0.82}
\definecolor{darkred}{RGB}{139,0,0}
\newcommand{\runin}[1]{\par\noindent\textbf{#1}\enspace}

\newcommand{\Top}{\ensuremath{T_{\mathrm{op}}}}
\newcommand{\Tev}{\ensuremath{T_{\mathrm{ev}}}}
\newcommand{\Thor}{\ensuremath{T_{\mathrm{h}}}}
\providecommand{\BibTeX}{
  {\rmfamily B\kern-.05em{\scshape i\kern-.025em b}\kern-.08em
  T\kern-.1667em\lower.7ex\hbox{E}\kern-.125emX}
}
\usepackage{xurl}
\usepackage{varioref}
\usepackage[hidelinks]{hyperref}
\usepackage[nameinlink,noabbrev]{cleveref}
\crefname{algorithm}{Algorithm}{Algorithms}
\Crefname{algorithm}{Algorithm}{Algorithms}
\begin{document}
\floatstyle{ruled}
\newfloat{algorithm}{tbp}{loa}
\floatname{algorithm}{Algorithm~}

\title{Counterfactual Screening of Rare Secondary Overloads in Power Grids\\
with a Stochastic Surrogate for Electromagnetic Transient Simulations}

\author{Ayrton Almada}
\affiliation{Program in Applied Mathematics and Department of Mathematics,
University of Arizona, Tucson, AZ 85721, USA}

\author{Laurent Pagnier}
\affiliation{Program in Applied Mathematics and Department of Mathematics,
University of Arizona, Tucson, AZ 85721, USA}

\author{Igal Goldshtein}
\affiliation{NOGA, Power System Operator of Israel, Haifa, Israel}

\author{Michael (Misha) Chertkov}
\affiliation{Program in Applied Mathematics and Department of Mathematics,
University of Arizona, Tucson, AZ 85721, USA}

\begin{abstract}
Power-system operators continually ask what would happen if a
component failed now. Single-phase-to-ground faults cleared by
single-pole tripping and automatic reclosure usually have no lasting
consequences. In rare cases, however, the transient excited during the
open-phase interval drives distant lines or transformers into overload
that persists after the faulted line is restored. We develop a
counterfactual screening framework for these secondary overloads,
combining a linear stochastic swing surrogate with cross-entropy
importance sampling. The surrogate represents the open-phase corridor
through an effective coupling factor and localized multiplicative
noise, with both parameters learned from electromagnetic-transient
(EMT) simulations rather than prescribed. Closed-form propagation of
the mean and covariance permits deterministic calibration without
Monte Carlo sampling. On the reduced WECC 240-bus system, a globally
calibrated model substantially improves finite-time agreement with
the EMT electromechanical state. The global models are calibrated at
reclosure and tested on disjoint faulted lines under a longer open-phase
interval. Line- and operating-state-specific models use additional
reclosure data from the assessment lines and provide smaller further
gains at the cost of more separate fits. The fitted coupling is lower than
the phase-counting baseline, consistent with the direction anticipated
by sequence-network analysis. Cross-entropy importance sampling
directs computational effort toward contingencies associated with
persistent overload and yields rare-event probability estimates
together with a map linking initiating faults to the elements they
endanger. By separating offline EMT calibration from repeated
surrogate evaluation, the method supports probabilistic contingency
screening and identifies where detailed high-fidelity assessment
is most valuable.
\end{abstract}


\maketitle

\section{Introduction}
\label{sec:intro}
\begin{figure*}[t]
\centering
\includegraphics[width=\textwidth]{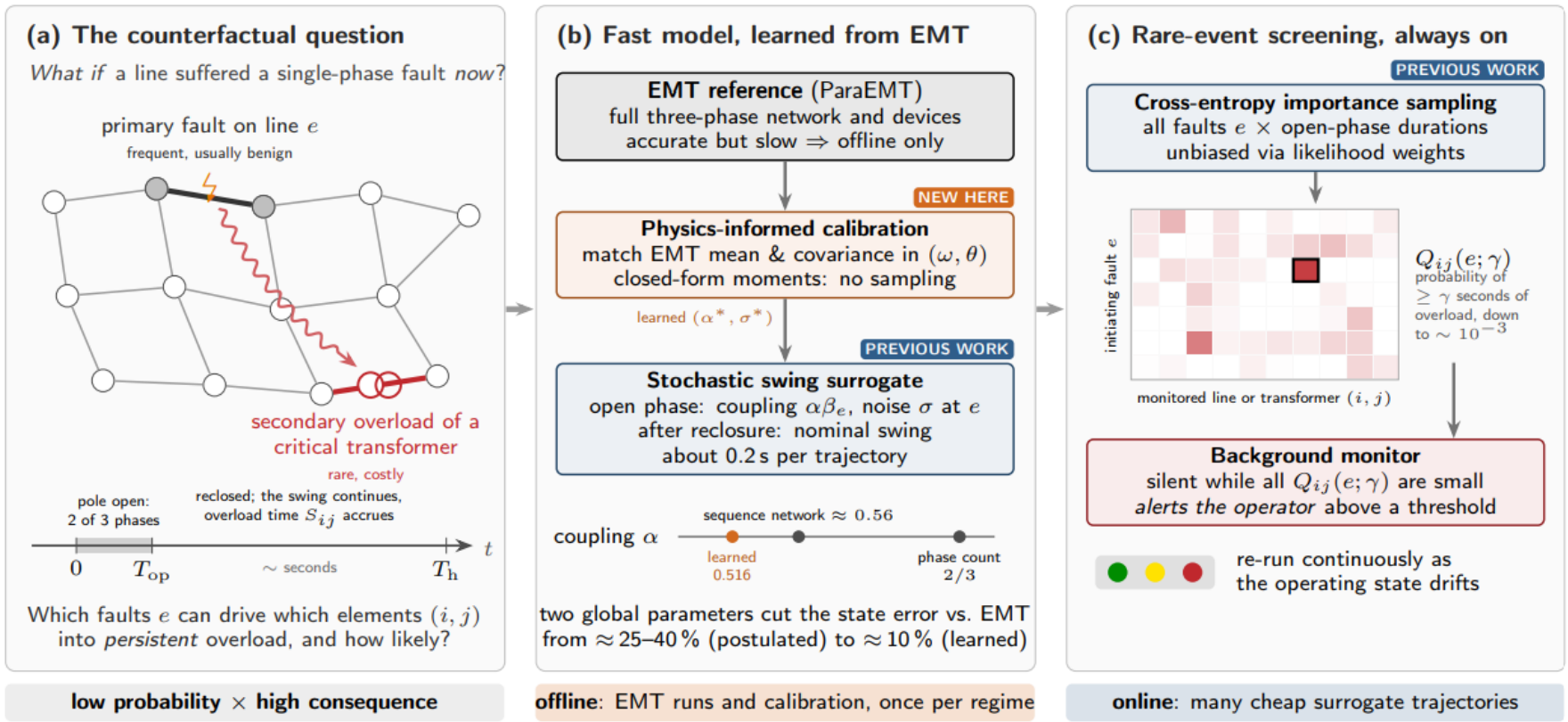}
\captionof{figure}{Logic of the paper.
(a)~The counterfactual question. A frequent and usually benign
single-phase fault on line $e$ is cleared by single-pole tripping and
reclosed after the open-phase interval $[0,T_{\rm op}]$. The
electromechanical swing it excites can drive a remote critical element
into persistent overload within the screening horizon $T_{\rm h}$.
(b)~The fast physics model. The open-phase coupling $\alpha$ and the
noise amplitude $\sigma$ of a stochastic swing surrogate are learned
offline from EMT (ParaEMT) trajectories by closed-form moment matching.
The learned coupling is compared with phase-counting and sequence-network
estimates.
(c)~Online, cross-entropy importance sampling over all initiating faults
turns surrogate trajectories into the contingency-to-vulnerable-element
map $Q_{ij}(e;\gamma)$. The intended tool runs in the background and
alerts the operator only when an entry exceeds a threshold.
Blue labels mark elements carried over from Ref.~\cite{almada2026realtime},
orange marks what is new here, and red marks secondary overload and the
operator alert. The network and map are schematic; quantitative results
are in Sec.~\ref{sec:results}.}
\label{fig:overview}
\end{figure*}

Most of the time, a transmission-system operator is watching a grid
in which nothing alarming is happening. The job is nevertheless
anticipatory: to know, before anything happens, which plausible events
would put the system in danger. Steady-state $N\!-\!1$ security assessment
answers this question for post-contingency equilibria and runs
continuously in control rooms \cite{ENTSOE_Operational_Handbook,NERC_Reliability_Standards}. It does not
answer it for the transient that connects the pre-fault and post-fault
states, even though that transient can drive other lines and transformers
above their limits for seconds. This paper develops a tool for that
missing, dynamic part of the question. The analysis is
\emph{counterfactual}. It is not a response to an ongoing emergency but an
estimate, made in advance and refreshed as conditions change, of what
\emph{would} happen if a given fault occurred now.

Guided by the operating experience of the Israeli transmission system
operator (NOGA), we
focus on the most common and least threatening transmission fault: a
single-phase-to-ground fault on an extra-high-voltage line, cleared by
single-pole tripping and automatic reclosure
\cite{montanari2009adaptative,godoy2012singlepole}. Opening only the
faulted pole keeps the two healthy phases in service. After a dead time of
a fraction of a second the pole recloses and, in the overwhelming majority
of cases, the event passes without consequence. Occasionally, however, the
electromechanical swing excited during the open-phase interval travels
through the network and pushes a remote critical element, typically a
transformer, into an overload that persists after the faulted line has
been restored. We call such a consequence a \emph{secondary contingency}.
It has three features that shape everything that follows. It is triggered
by a frequent, benign primary event. It is rare. And when it happens it is
costly, since a sustained transformer overload can lead to damage,
protective tripping, or the first step of a cascade.

Screening for secondary contingencies is therefore a rare-event
problem posed under real-time constraints, and three requirements pull
against each other. The first is \emph{fidelity}. The open-phase interval
is an unbalanced three-phase transient, and electromagnetic-transient
(EMT) simulation is the industry-accepted description of it
\cite{van2007voltage,xiong2024open}. Even a parallel EMT code needs about
half a minute of computation per simulated second on the reduced 240-bus
model of the Western Electricity Coordinating Council (WECC) system
\cite{xiong2023paraemt}. The second is \emph{rarity}: resolving a
probability of order $10^{-3}$ to within ten percent by direct Monte Carlo
requires of order $10^{5}$ trajectories. The third is \emph{extent and
latency}: every line is a possible initiating fault, and the whole
assessment must be repeated as load, dispatch, and renewable output drift.
Direct EMT Monte Carlo meets the first requirement and misses the other two
by orders of magnitude.

In previous work \cite{almada2026realtime} we addressed rarity and
extent. Following the swing-equation tradition of transient-stability
analysis \cite{wood2013power,sauer2017power,pavella2012transient}, we
represented the post-fault dynamics by a linear stochastic swing model in
which the faulted corridor has a reduced coupling during the open-phase
interval and carries multiplicative noise localized at its endpoints.
Combined with cross-entropy (CE) importance sampling
\cite{rubinstein_optimization_1997,de2005tutorial,owen2019importance},
this model screens many initiating single-phase faults in parallel and
concentrates computational effort on the rare scenarios that produce long
overloads. The weak point of that construction was its physics. The
reduced coupling and the noise amplitude were postulated, not derived from
or tested against a model that operators trust. Whether such a light
model can stand in for the EMT description of the open-phase transient was
left open.

This paper identifies the surrogate's effective parameters from
ParaEMT reference trajectories on the reduced WECC 240-bus system
\cite{xiong2024open}. Calibration matches the electromechanical state
at reclosure; subsequent assessment tests the globally fitted model on
disjoint faulted lines under a longer open-phase interval. A hierarchy
of global, line-specific, and line--operating-state-specific fits examines
the value of additional calibration information. The comparison uses the
phase-angle and frequency variables retained by the swing model, with
parameter identification based on its deterministic moment equations.

The calibrated surrogate is then used to estimate overload probabilities
and construct a \emph{contingency-to-vulnerable-element map}: for each
initiating fault, the probability that a monitored line or transformer
accumulates more than a prescribed time in overload. This is the
operator-facing output of the proposed background screening workflow in
Fig.~\ref{fig:overview}.

Section~\ref{sec:methods} develops the model, calibration, and probability
estimation; Sec.~\ref{sec:wecc} specifies the experiments.
Sections~\ref{sec:results} and~\ref{sec:conclusion} present the findings
and their interpretation, with scope and limitations collected in
Sec.~\ref{subsec:scope_limitations}.

\section{Methods}\label{sec:methods}

\begin{table*}[t]
\caption{Principal notation and time conventions. Physical calibration
parameters, scenario-law parameters, and sampling diagnostics have distinct
roles. The same time symbol may take different values in different experiments.}
\label{tab:global_notation}
\label{tab:terminology}
\small
\renewcommand{\arraystretch}{1.13}
\begin{tabular}{@{}p{0.28\textwidth}p{0.69\textwidth}@{}}
\toprule
\textbf{Symbol} & \textbf{Meaning}\\
\midrule
$\mathcal G=(\mathcal V,\mathcal E)$; $n_{\rm E}=|\mathcal E|$ & Network graph, bus and line sets, and number of initiating lines; $n_{\rm E}=329$ in screening.\\
$e=(i,j)$; $a\in\mathcal N_0(e)$ & Faulted line and a noise-bearing endpoint; $\mathcal N_0(e)=\{i,j\}$.\\
$z$; $\chi$ & Pre-fault operating state and its dimensionless uniform injection multiplier.\\
$X_t=(\omega_t,\theta_t)^{\mathsf T}$ & Reduced state: frequency deviations and phase angles in the fixed reference.\\
$\mathbf M$; $\mathbf e_a$ & Inertia matrix and the canonical vector for endpoint $a$.\\
$A$, $A_e(\alpha)$; $A_{e,\mathrm I}$ & Nominal drift, open-phase drift, and equivalent It\^o drift.\\
$G_{e,a}(\sigma)$; $W_t^a$ & Endpoint-local diffusion matrix and its independent Wiener driver.\\
$\alpha$ & Open-phase coupling factor: $\beta_e\mapsto\alpha\beta_e$ during $[0,\Top]$ only.\\
$\sigma$ & Common amplitude of endpoint-localized multiplicative noise; $[\sigma]=\mathrm{s}^{-1/2}$.\\
$\Theta$ & Fitted parameter collection.\\
$\mathrm M_0,\ldots,\mathrm M_4$ & Calibration-hierarchy model names, not inertia matrices.\\
$t=0$ & Instant the faulted pole opens; fault clearing is treated as instantaneous.\\
$\Top$ & Open-phase duration (auto-reclosing dead time).\\
$\Tev$ & State-evaluation time: $0.8\,\mathrm{s}$ at reclosure for every fitted hierarchy level; $2.5\,\mathrm{s}$ after a $1.5$-s open-phase interval for their common assessment.\\
$\Thor$ & Upper limit of the overload integrals, fixed at $10\,\mathrm{s}$ in screening; distinct from the calibration target time.\\
$H$; $y_{e,z}^{\rm EMT}$; $\mathcal T_{e,z}^{\rm EMT}$ & Observation convention, projected EMT state, and finite-time reference state-transfer map.\\
$\mu_\Theta$, $\Sigma_\Theta$ & Surrogate mean and state covariance.\\
$W_{e,H}^{(q)}$; $\iota_{e,a}^{(q)}$ & Observation-space neighborhood mask and its binary bus-selection entries.\\
$\widehat\Sigma_e^{\rm EMT}$ & Covariance target formed across calibration operating states at fixed line.\\
$\Sigma_{\Theta,e,z}^{(0)}$, $\widetilde\Sigma_{\Theta,e,z}^{(0)}$ & Endpoint-masked predicted covariance and its diagonally regularized version.\\
$J_{\rm map}^{(q)}$; $\mathcal R$; $\zeta_{\rm cov},\zeta_{\rm reg}$ & Calibration objective, scalar regularizer, and objective weights.\\
$C$; $\lambda\in\operatorname{spec}(K_{e,\Theta})$ & Matrix argument of the second-moment operator and an eigenvalue of its matrix representation.\\
$\rho=X_0/X_1$ & Zero- to positive-sequence reactance ratio in Appendix~\ref{app:sequence_network}.\\
\midrule
$S_{ij}$; $S$; $\gamma\geq0$ & Accumulated line exposure, summed network exposure, and the prescribed strict-exceedance threshold.\\
$Q(\gamma)$, $Q_{ij}(\gamma)$ & Global and line-specific probabilities under the nominal scenario law.\\
$Q(e;\gamma)$, $Q_{ij}(e;\gamma)$ & Corresponding probabilities conditional on initiating fault $e$.\\
$\tau$; $\tau_0$; $r$ & Sampled open-phase duration, nominal exponential mean, and shared proposal mean, all in seconds.\\
$\pi_e$; $\phi_e$; $\nu$ & Nominal line probability, proposal line probability, and proposal-parameter vector.\\
$N_{\rm CE}$, $K$, $N$ & Samples per adaptation iteration, number of iterations, and fresh final-batch size.\\
$\rho_{\rm CE}$; $\kappa$; $n_{\rm elite}$ & CE elite fraction, previous-parameter retention factor, and elite count.\\
$\ell$; $k$; $\mathcal I_\ell$ & Adaptation iteration, sample index, and selected elite index set.\\
$b_k^{(\ell)}$; $s_\ell^{\rm elite}$ & Binary elite-membership indicator and smallest selected elite score.\\
$w_k$; $\widetilde w_k$; $Y_k$ & Raw likelihood ratio, normalized diagnostic weight, and weighted event contribution.\\
$\widehat Q$; $\widehat{\operatorname{se}}$; $\operatorname{RelSE}$ & Probability estimate, estimated standard error, and relative sampling standard error.\\
$\mathrm{RSE}_{\rm map}$ & Relative finite-time state error; distinct from sampling uncertainty.\\
$\mathbb P$, $\mathbb E$, $\exp$ & Probability, expectation, and the exponential function.\\
\bottomrule
\end{tabular}
\end{table*}

Table~\ref{tab:global_notation} collects the notation and time conventions
used in the main text and appendices. The method separates representation from estimation: EMT simulations
calibrate the reduced dynamics, and importance sampling estimates overload
probabilities under those fixed dynamics. Calibration, state assessment,
and screening use distinct time settings. The state variables, switching
rule, and observation convention remain the same; the fitted physical
parameters are held fixed when the screening sampler varies fault
location and open-phase duration.

\subsection{Contingency and Reduced Dynamics}
\label{sec:surrogate}
\label{subsec:terminology}
\label{subsec:postclearing}
The modeled contingency is a single-phase-to-ground fault cleared by
single-pole tripping and followed by successful automatic reclosure.
The reduced description starts at pole opening and retains the open-phase
and post-reclosure electromechanical response.

\runin{Fault inception and clearing.} Protection opens the faulted
pole of line $e=(i,j)$. Fault inception and breaker operation occur over
tens of milliseconds, below the electromechanical time scale retained by
the surrogate. Clearing is therefore treated as instantaneous, and
$t=0$ denotes the instant of pole opening. The reduction begins from the
corresponding electromechanical state and does not resolve the preceding
fault or breaker transient.

\runin{Open-phase interval $[0,\Top]$.} One phase is de-energized,
but the other two remain in service, so the line continues to transmit
power. The interval permits de-ionization of the secondary arc sustained
by capacitive and inductive coupling from the energized phases; premature
reclosure can restrike the fault. In protection terminology, $\Top$ is
the auto-reclosing \emph{dead time} \cite{IEC60050_441}. We use
\emph{open-phase duration} to distinguish it from the evaluation and
screening times. During this interval, the effective positive-sequence
coupling changes from $\beta_e$ to $\alpha\beta_e$, with $\alpha<1$.

\runin{Reclosure and post-reclosure interval $(\Top,\Thor]$.}
At $t=\Top$, the open pole is reclosed and the nominal coupling
$\beta_e$ is restored \cite{montanari2009adaptative}. The machines do
not instantaneously return to equilibrium: their ensuing swing can cause
transient overloads on other elements. The transferred electromechanical
state therefore matters even after the network configuration is restored.

\runin{Evaluation time $\Tev$ and screening horizon $\Thor$.}
The surrogate state is compared with its EMT reference at the
experiment-specific time $\Tev$; overload duration is accumulated up to
$\Thor$. State assessment at or after reclosure has
$0<\Top\leq\Tev$. Calibration uses the boundary case
$\Top=\Tev=0.8\,\mathrm{s}$, so the target is the state attained at
reclosure. The main hierarchy comparison instead uses
$\Top=1.5\,\mathrm{s}$ and $\Tev=2.5\,\mathrm{s}$, one second after
reclosure. Screening integrates exposure through $\Thor=10\,\mathrm{s}$
and does not fix every sampled duration at a calibration or test value.
A full-trajectory illustration may also end at $10\,\mathrm{s}$ without
making that endpoint the calibration time. Varying $\Tev$ changes the
observation time; varying $\Top$ changes the switching event.
Table~\ref{tab:global_notation} records these roles. Throughout, fault clearing
initiates the open-phase interval at $t=0$; $t=\Top$ denotes reclosure,
not a second fault-clearing event.

\runin{Why a stochastic reduced model?}
The swing variables retain the phase and frequency response used by the
screening calculation, but not the three-phase electrical and device
states resolved by ParaEMT. Their influence on the retained variables is represented by an
effective coupling and a linear stochastic correction, whose parameters
are identified from the EMT response.

The two fitted terms address different features of the reduction.
The coupling factor represents the change in transfer through the
open-phase corridor. The diffusion matrices act only at its endpoints,
so the model introduces the correction at the location of the disturbance
and retains the network dynamics that transmit its effects elsewhere.
Linearity in the state gives the closed moment equations used for
calibration.

The stochastic correction represents unresolved variability within the
prescribed operating-state family. Its empirical assessment combines the
open-phase comparison in Fig.~\ref{fig:open_phase_emt_surrogate} with
the deterministic--stochastic comparison in the calibration hierarchy.
Section~\ref{subsec:moment_calibration} specifies the across-operating-state
covariance target used to identify its amplitude.

For a faulted line $e=(i,j)$ and a pre-fault operating state $z$,
the surrogate evolves the electromechanical state through the
open-phase and post-reclosure intervals of
Sec.~\ref{subsec:terminology}.

The reduced electromechanical state is
\begin{equation}
    X_t =
    \begin{bmatrix}
        \omega_t\\
        \theta_t
    \end{bmatrix},
\end{equation}
where $\omega_t$ and $\theta_t$ are the frequency-deviation and phase
vectors over the bus set $\mathcal V$. The formulation preserves the
reduced WECC network graph and transmission-line connectivity rather than
Kron-reducing the system to generator internal nodes. Thus, the graph
$\mathcal G$ used for the neighborhoods in
\eqref{eq:line_distance}--\eqref{eq:local_neighborhood} is the graph
used to construct the swing dynamics. Phases are absolute coordinates
relative to a fixed reference bus, not deviations from the pre-fault
phase vector; $X_0(z)$ denotes the equilibrium state for operating
condition $z$ in that reference convention.

The symbols $\alpha$ and $\sigma$ denote the coupling and common
endpoint-noise amplitude used for a specified scenario. Their sharing
across scenarios is determined by the calibration hierarchy in
Sec.~\ref{subsec:calibration_protocol}, under \emph{Comparison hierarchy},
and Fig.~\ref{fig:sequence_control_design}. They are global constants in
$\mathrm M_1/\mathrm M_2$, line-specific values $(\alpha_e,\sigma_e)$ in
$\mathrm M_3$, and line--operating-state-specific values
$(\alpha_{e,z},\sigma_{e,z})$ in $\mathrm M_4$. The baseline
$\mathrm M_0$ fixes $(\alpha,\sigma)=(2/3,0)$, and $\mathrm M_1$ has
$\sigma=0$. In every case, one amplitude is shared by the two endpoints
of a given faulted line. We suppress these calibration indices in the
dynamical equations; the scenario indices on $A_e$ and $G_{e,a}$ still
identify the modified corridor and the endpoint where noise acts.

Across the open-phase and post-reclosure intervals, the dynamics are
\begingroup
\footnotesize\let\small\relax\thinmuskip=2mu\medmuskip=2mu\thickmuskip=2mu
\begin{equation}
\small
\begin{cases}
    dX_t
    =
    \left(
        A_e(\alpha)X_t + P(z)
    \right)dt
    +
    \sum_{a\in\mathcal N_0(e)}
    G_{e,a}(\sigma)X_t\circ dW_t^a:\\
    0\leq t\leq \Top,\\
    \dot{X_t}=AX_t + P(z):
    \Top\leq t\leq \Tev,
\end{cases}
\label{eq:surrogate_sde}
\end{equation}
\endgroup
where
\begin{equation}
    \mathcal N_0(e) = \{i,j\}
    \label{eq:fault_endpoints}
\end{equation}
is the endpoint set of the faulted line, equivalently its zero-radius
graph neighborhood. During $[0,\Top]$, the
multiplicative forcing represents variability localized to the affected
corridor. The equation is interpreted in the Stratonovich sense, as
indicated by $\circ dW_t^a$, and the $W_t^a$ are independent Wiener
processes. The vector $P(z)$ denotes the power-injection forcing for the
pre-fault operating state. At reclosure, the stochastic forcing is removed
and the nominal drift $A$ is restored; the attained electromechanical
state supplies the initial condition for the post-reclosure interval.
The linear description therefore yields a finite-time map composed across
the switching event, rather than one time-invariant faulted-network
propagator.

The contingency-dependent matrix $A_e(\alpha)$ is constructed
from the swing dynamics by modifying the coupling of the affected line,
\begin{equation}
    \beta^{\mathrm{op}}_e
    =
    \alpha\beta_e,
    \label{eq:effective_coupling}
\end{equation}
where $\beta_e$ is its nominal coupling. The modification applies only
during the open-phase interval. Reclosure restores $\beta_e$ and
$A=A_e(1)$, so the calibrated parameters refer to the same piecewise
network configuration used in trajectory evaluation.

The parameter $\alpha$ is the effective positive-sequence coupling
factor of the open-phase corridor. Counting the two intact phases gives
the physically transparent baseline $\alpha=\frac{2}{3}$, but not an
exact sequence-network result. Negative- and zero-sequence contributions
modify the positive-sequence transfer under a one-conductor-open condition,
providing a physical basis for a coupling below the phase-counting value.
The value used by the surrogate is identified from the electromechanical
response rather than imposed by phase counting; its comparison with the
sequence-network estimate is reported in
Sec.~\ref{subsec:global_wecc_calibration}.

The noise acts on the phase coordinates of the two endpoints. Let
$\mathbf{e}_a$ be the $a$th canonical vector in
$\mathbb R^{|\mathcal V|}$. For $a\in\mathcal N_0(e)$,
\begin{equation}
    G_{e,a}(\sigma)
    =
    \begin{bmatrix}
        0 & 0\\
        0 & \sigma\mathbf M^{-1}\mathbf{e}_a\mathbf{e}_a^{\mathsf T}
    \end{bmatrix}.
    \label{eq:noise_matrix}
\end{equation}
Here $\mathbf M$ is the inertia matrix of the retained swing model.
Each diffusion matrix has a single endpoint contribution in the phase block. Since
$dW_t$ has units $\mathrm{s}^{1/2}$, the amplitude has units
$[\sigma]=\mathrm{s}^{-1/2}$; this is not a dimensionless perturbation. Because the multiplicative term uses phase coordinates,
its specification is tied to the fixed phase reference introduced above.
The same reference must be retained in the EMT projection, initial states,
and surrogate evaluation; common-mode cancellation in the overload
observable does not make the stochastic state equation reference-invariant.

The two parameters are identified jointly. The coupling $\alpha$
controls the deterministic open-phase interaction, while $\sigma$
controls localized variability and contributes to the equivalent It\^o
drift through $\frac12\sum_{a\in\mathcal N_0(e)}G_{e,a}^2$. The baseline
$\alpha=\frac23$ provides a warm start, not a constraint on the fit.
The stochastic correction can consequently influence both the transferred
mean state and its covariance, as made explicit in
Sec.~\ref{subsec:moment_calibration}.

The sequence-network construction in Appendix~\ref{app:sequence_network}
provides a physical motivation for the coupling correction. The hierarchy
and the status of the sequence-informed control are specified in
Sec.~\ref{subsec:calibration_protocol}.

\subsection{EMT Reference and Moment-Based Calibration}
\label{sec:learning}
\label{sec:emt_reference}
\label{subsec:emt_projection}
\label{subsec:moment_calibration}
\label{subsec:map_error}
\runin{Shared observation space.}
ParaEMT provides the reference trajectories from which the reduced
finite-time map is identified \cite{xiong2024open}. Its three-phase
network and device model supplies the high-fidelity response, which is
projected onto the observation space shared with the swing surrogate.
Appendix~\ref{app:reference_details} identifies the reference
implementation and case inputs.

Let $e=(i,j)$ denote the faulted line and $z$ the pre-fault operating
state. ParaEMT generates a high-dimensional trajectory
\begin{equation}
    x^{\mathrm{EMT}}_{e,z}(t),
    \qquad
    0\leq t\leq \Tev ,
\end{equation}
where $\Tev$ is the evaluation time defined in
Sec.~\ref{subsec:terminology}. Its electromechanical observation is
\begin{equation}
    y^{\mathrm{EMT}}_{e,z}(t)
    =
    Hx^{\mathrm{EMT}}_{e,z}(t),
    \qquad
    0\leq t\leq \Tev ,
    \label{eq:emt_observation}
\end{equation}
where $H$ extracts the phase and frequency variables represented by the
surrogate. The same observation convention is used when evaluating the
reduced state, so the mismatch concerns corresponding physical variables.

The projection retains the electromechanical response through which the
reference model's fast electromagnetic and device dynamics influence the
calibration target. The selected observations determine the dimension of
the comparison.

The finite-time reference map is
\begin{equation}
    \mathcal{T}^{\mathrm{EMT}}_{e,z}
    :
    X_0(z)
    \longmapsto
    y^{\mathrm{EMT}}_{e,z}(\Tev ),
    \label{eq:emt_map}
\end{equation}
where $X_0(z)$ is the reduced electromechanical state associated with the
pre-fault equilibrium. Fixing the initial state, switching convention, and
observation time makes the map a direct comparison of the reference and
reduced responses to the same scenario. Its approximation is the
calibration objective. In the reported fit this target is evaluated at
$\Tev=\Top=0.8\,\mathrm{s}$; the subsequent assessment propagates the
fitted model to its separately specified observation time.

\runin{Deterministic moment propagation.}
Linearity with multiplicative noise gives a closed deterministic system
for the first two moments after conversion to the equivalent It\^o form.
For parameters $\Theta$, define the mean
\begin{equation}
    \mu_{\Theta}(t)
    =
    \mathbb{E}[X_t],
\end{equation}
and covariance
\begin{equation}
    \Sigma_{\Theta}(t)
    =
    \mathbb{E}
    \left[
        (X_t-\mu_{\Theta}(t))
        (X_t-\mu_{\Theta}(t))^{\mathsf T}
    \right].
\end{equation}
The equivalent It\^o drift during the open-phase interval is
\begin{equation}
    A_{e,\mathrm I}(\alpha,\sigma)
    =
    A_e(\alpha)
    +
    \frac{1}{2}
    \sum_{a\in\mathcal N_0(e)}
    G_{e,a}^2(\sigma).
\end{equation}
The second term is the noise-induced drift associated with the
Stratonovich convention
\cite{gardiner2009stochastic,oksendal2003stochastic,mannella2012ito}.
Consequently, changing $\sigma$ affects the mean as well as the
covariance; its role is not restricted to adding fluctuation around an
unchanged deterministic trajectory. The mean satisfies
\begingroup
\small
\begin{equation}
    \dot{\mu}_{\Theta}(t)
    =
    \begin{cases}
        \left[
            A_e(\alpha)
            +
            \dfrac{1}{2}
            \displaystyle\sum_{a\in\mathcal N_0(e)}
            G_{e,a}^2(\sigma)
        \right]\mu_{\Theta}(t)
        +P(z):\\
        0\leq t\leq \Top,\\[8pt]
        A\mu_{\Theta}(t)+P(z):\Top\leq t\leq \Tev.
    \end{cases}
    \label{eq:mean_dynamics}
\end{equation}
\endgroup
and the covariance satisfies
\begingroup
\footnotesize
\begingroup
\footnotesize\let\small\relax
\begin{equation}
\small
    \dot{\Sigma}_{\Theta}(t)
    =
    \begin{cases}
        A_{e,\mathrm I}\Sigma_{\Theta}
        +
        \Sigma_{\Theta}A_{e,\mathrm I}^{\mathsf T}
        +
        \displaystyle\sum_{\mathclap{a\in\mathcal N_0(e)}}
        G_{e,a}
        \left(
            \Sigma_{\Theta}
            +
            \mu_{\Theta}\mu_{\Theta}^{\mathsf T}
        \right)
        G_{e,a}^{\mathsf T}:\\
        0\leq t\leq \Top,\\[8pt]
        A\Sigma_{\Theta}
        +
        \Sigma_{\Theta}A^{\mathsf T}:\Top\leq t\leq \Tev,
    \end{cases}
    \label{eq:covariance_dynamics}
\end{equation}
\endgroup
\endgroup
The moments at reclosure initialize their deterministic continuation
under $A$. In calibration, evaluation stops at reclosure; in the
post-reclosure tests, the moments are continued to the specified $\Tev$.
This retains the switching event explicitly and avoids Monte Carlo
sampling when evaluating candidate parameters
\cite{gardiner2009stochastic,oksendal2003stochastic}.

The effective coupling and noise-induced drift jointly determine the
mean-based metric in \eqref{eq:rse_map}. This supplies the mechanism
by which $\mathrm M_2$ can improve the finite-time map beyond the deterministic
fit $\mathrm M_1$: $\alpha$ adjusts the corridor interaction, while the
Stratonovich correction permits $\sigma$ to modify the transferred
mean. Their joint interpretation allows part of the effective correction
to be represented by noise-induced drift rather than by $\alpha$ alone.
Appendix~\ref{app:ms_stability} records the second-moment operator and
the spectral diagnostic for the reduced stochastic dynamics.

\runin{Endpoint-local objective.}
The primary global fit uses only the two endpoints of the faulted line.
The mask $W_{e,H}^{(0)}$ selects their observed phase and frequency
components, and $\widehat\Sigma_e^{\rm EMT}$ denotes the scatter of the
projected EMT state across the training operating-state family at fixed
fault $e$, evaluated at $\Tev=\Top=0.8\,\mathrm{s}$ for the training
runs. The joint objective is
\begingroup
\footnotesize
\begingroup
\footnotesize\let\small\relax
\begin{equation}
\small
\begin{aligned}
& J_{\rm map}^{(q)}(\Theta)
\doteq 
\sum_{(e,z)\in{\cal D}_{\rm train}}
\left\|
W_{e,H}^{(q)}
\left(
H\mu_\Theta(\Tev ;e,z)
-
y^{\rm EMT}_{e,z}(\Tev )
\right)
\right\|_2^2
\\
&
+
\zeta_{\rm cov}
\sum_{(e,z)\in{\cal D}_{\rm train}}
\left\|
W_{e,H}^{(q)}
\left(
H\Sigma_\Theta(\Tev ;e,z)H^{T}
-
\widehat\Sigma^{\rm EMT}_{e}
\right)
W_{e,H}^{(q)}
\right\|_{F}^{2}\\ & 
+
\zeta_{\rm reg}\mathcal R(\Theta).
\end{aligned}
\label{eq:calibration_objective}
\end{equation}
\endgroup
\endgroup
The first term matches the mean transferred state, the second matches
its covariance target, and $\mathcal R(\Theta)$ regularizes the predicted
covariance. The weights are $\zeta_{\rm cov}=1.0$ and
$\zeta_{\rm reg}=0.5$; the absolute log-determinant regularizer uses
the diagonal shift $\epsilon=0.01$. Their explicit definitions and
mask construction are given in Appendix~\ref{subsec:graph_local}.

ParaEMT supplies one deterministic trajectory for each $(e,z)$, and
$\widehat\Sigma_e^{\rm EMT}$ is computed across the training operating
states. The fitted $\sigma$ therefore represents operating-point
variability unresolved by the global reduced description. Only training
states from the global calibration-line set enter this target; observations
from the global testing lines are not used to identify the shared
parameters of $\mathrm M_1$ or $\mathrm M_2$. Models $\mathrm M_3$ and $\mathrm M_4$ instead receive
additional calibration observations from those lines at the same
$0.8$-s reclosure target, as specified in
Sec.~\ref{subsec:calibration_protocol}. The later observations used to
score Fig.~\ref{fig:boxplots} are not fitting targets for any model.
The restriction $q=0$ applies to parameter identification, not to the
propagation of the disturbance: network-wide responses remain available
for validation and screening.

\runin{Validation metric.}
The principal validation metric compares the surrogate mean and the
projected EMT state at the same evaluation time,
\begin{equation}
    \mathrm{RSE}_{\mathrm{map}}(e,z)
    =
    \frac{
        \left\|
            H\mu_{\Theta}(\Tev ;e,z)
            -
            y^{\mathrm{EMT}}_{e,z}(\Tev )
        \right\|_2
    }{
        \left\|
            y^{\mathrm{EMT}}_{e,z}(\Tev )
        \right\|_2
    }.
    \label{eq:rse_map}
\end{equation}
The numerator measures discrepancy in the observed electromechanical
state; the denominator expresses that discrepancy relative to the
reference state's norm. Both the phase-angle and frequency components
selected by $H$ enter this metric; it is not a phase-only error unless
$H$ selects only phases. Within a given comparison, all hierarchy levels
use the same variables, reference convention, switching duration, and
observation time. This common scoring protocol does not imply identical
fitting data: the global and scenario-specific models use the allocations
in Sec.~\ref{subsec:calibration_protocol}. All fitted coefficients are
held fixed when the $0.8$-s calibration setting is replaced by the
longer-duration assessment setting.

In Fig.~\ref{fig:boxplots}, the finite-time state metric is evaluated at
$\Tev=2.5\,\mathrm{s}$ after reclosure at $\Top=1.5\,\mathrm{s}$,
with the coefficients fitted at the $0.8$-s calibration target held fixed.
Noise affects $H\mu_\Theta$ through the mean dynamics in
Sec.~\ref{subsec:moment_calibration}. Covariance information enters
parameter identification through \eqref{eq:calibration_objective};
validation then scores the projected mean state.

\subsection{Overload Observables and Conditional Risks}
\label{sec:probabilistic_screening}
\label{subsec:dynamic_overload}
\label{subsec:rare_event_probability}
Let $\theta(t)$ be the phase-angle vector of a sampled surrogate trajectory
from \eqref{eq:surrogate_sde}. For a monitored line
$\{i,j\}\in\mathcal{E}_m$, the linearized transfer
$\beta_{ij}[\theta_i(t)-\theta_j(t)]$ is compared with its allowable
power-flow threshold $\bar p_{ij}$. The line-specific overload indicator is
\begin{equation}
    S_{ij}(\theta_{0\rightarrow \Thor})
    =
    \int_0^{\Thor}
    \mathbb{I}
    \left\{
        \left|
            \beta_{ij}
            \left[
                \theta_i(t)-\theta_j(t)
            \right]
        \right|
        >
        \bar{p}_{ij}
    \right\}
    dt,
    \label{eq:line_overload_indicator}
\end{equation}
where $\mathbb I\{\cdot\}$ equals one when the stated limit is exceeded and
zero otherwise. Thus, $S_{ij}$ is a cumulative duration: separated
periods of overload contribute to the same integral. The monitored element
need not be the faulted line. This distinction allows a local initiating
contingency to be assessed through the persistent stress it produces
elsewhere in the network.

The system-level indicator sums exposure over the monitored set,
\begin{equation}
    S(\theta_{0\rightarrow \Thor})
    =
    \sum_{\{i,j\}\in\mathcal{E}_m}
    S_{ij}(\theta_{0\rightarrow \Thor}),
    \label{eq:global_overload_indicator}
\end{equation}
where $\mathcal{E}_m\subseteq\mathcal{E}$. Simultaneous violations on
different elements contribute separately, so $S$ measures accumulated
line-overload time rather than the duration for which at least one line
is overloaded. It can consequently exceed $\Thor$. Neither indicator
weights a violation by the amount of its excess above the threshold;
their purpose is to distinguish persistent excursions from brief crossings
that may have comparable peak phase separations.

Equation~\eqref{eq:line_overload_indicator} uses a linearized-flow
screening approximation. Common angular displacements cancel from the
phase differences. Figure~\ref{Point3} and
Appendix~\ref{app:large_angle} examine the corridor separations that
determine its small-angle interpretation; the scope of the indicator is
set out in Sec.~\ref{subsec:scope_limitations}.

For a prescribed overload-duration threshold $\gamma\geq0$, define
\begin{equation}
    Q_{ij}(\gamma)
    =
    \mathbb{P}
    \left[
        S_{ij}(\theta_{0\rightarrow \Thor})
        >
        \gamma
    \right].
    \label{eq:line_probability}
\end{equation}
This is the probability that monitored line $\{i,j\}$ accumulates more
than $\gamma$ seconds above its flow limit. The corresponding global
probability is
\begin{equation}
    Q(\gamma)
    =
    \mathbb{P}
    \left[
        S(\theta_{0\rightarrow \Thor})
        >
        \gamma
    \right].
    \label{eq:global_probability}
\end{equation}
Here, the event concerns accumulated exposure over the monitored network,
not the union of line-specific events at the same threshold. In particular,
$Q(\gamma)$ is not obtained by adding the $Q_{ij}(\gamma)$. The strict
comparison also defines $Q(0)$ as the probability of any positive overload,
rather than the certain event associated with a nonnegative indicator.

The probability law specifies both the initiating scenario and the
surrogate uncertainty conditional on it. A scenario comprises a faulted
line $e$ and its open-phase duration $\Top$; the stochastic drivers
then determine a realization of the electromechanical response. The
nominal scenario law remains fixed when the sampling algorithm changes.
The conditional risks defined below answer a different question from
averaging over possible fault locations.

\runin{Conditional risk scores.}
For an initiating fault $e$ and threshold $\gamma$, define
\begin{equation}
    Q(e;\gamma)
    =
    \mathbb P\!\left[
        S\!\left(\theta_{0\rightarrow\Thor}^{(e)}\right)>\gamma
    \right],
    \label{eq:contingency_risk}
\end{equation}
The probability averages over the nominal open-phase duration and
surrogate drivers conditional on $e$. Its strict comparison follows
\eqref{eq:global_probability}. Sorting $Q(e;\gamma)$ gives a ranking
by the probability of exceeding a specified total line-overload time.
Because $S$ includes the monitored set, the score can include exposure
on the faulted element as well as secondary exposure elsewhere.

For an individual monitored line, define
\begin{equation}
    Q_{ij}(e;\gamma)
    =
    \mathbb P\!\left[
        S_{ij}\!\left(\theta_{0\rightarrow\Thor}^{(e)}\right)
        >\gamma
    \right],
    \qquad \{i,j\}\in\mathcal{E}_m,
    \label{eq:conditional_line_risk}
\end{equation}
The strict comparison is the convention in
\eqref{eq:line_probability} and \eqref{eq:conditional_risk_estimator}.
The collection
\begin{equation}
    \left\{
        Q_{ij}(e;\gamma)
    \right\}_{\{i,j\}\in\mathcal{E}_m}
    \label{eq:contingency_vulnerability_map}
\end{equation}
forms one row of the contingency-to-vulnerable-element map. Entries for
$\{i,j\}\neq e$ identify secondary exposure. Conversely, examining
a monitored element across initiating faults identifies whether it is
repeatedly vulnerable or associated mainly with a particular contingency.
The global score is evaluated from $S$, not from a sum of these
line-specific event probabilities.

Appendix~\ref{app:conditional_risk} gives the conditional estimator and
its normalization by the nominal fault probability. This normalization
separates a risk conditional on $e$ from the contribution of $e$ to an
ensemble that samples many possible fault locations. The learned CE masses
control sampling effort and are not themselves conditional overload
probabilities.

\subsection{Cross-Entropy Importance Sampling}
\label{subsec:importance_sampling}
\label{subsec:role_of_calibration}
As the overload threshold increases, most nominal samples may contribute
zero to the event indicator. Importance sampling reallocates the
trajectory budget without changing either the event or the calibrated
physical model. The stochastic driving law is the same under the nominal
and proposal simulations; only fault location and duration are resampled.

Direct Monte Carlo draws the faulted line and open-phase duration from
the nominal joint law $p_Z$. CE adaptation instead fits a proposal to
high-score contingencies and then freezes the proposal for a separate
estimation batch. The WECC implementation retains both the learned
categorical line probabilities and a shared duration parameter. With
$\tau_k$ denoting the sampled open-phase duration, the estimator is
\begin{equation}
    \widehat Q(\gamma)
    =
    \frac{1}{N}\sum_{k=1}^N
    \mathbb I\!\left\{
        S\!\left(\theta_{0\rightarrow\Thor}^{(k)}\right)>\gamma
    \right\}
    \frac{p_Z(e_k,\tau_k)}{q_{\nu^\ast}(e_k,\tau_k)}.
    \label{eq:cem_estimator}
\end{equation}
The ratio $p_Z/q_{\nu^\ast}$ corrects both line selection and duration
sampling. Positive proposal support on the nominal support, an unchanged
stochastic driving law, and a fresh estimation batch preserve the nominal
target. The joint law and likelihood ratio are specified in
Appendix~\ref{subsec:cross_entropy}
\cite{rubinstein_optimization_1997,de2005tutorial,owen2019importance}.
The elite updates and Algorithm~\ref{alg:n1plus} are given in the same
appendix. The duration parameter $r$ is an exponential mean, measured in
seconds; the nominal mean $\tau_0=0.1\,\mathrm{s}$ remains fixed during
proposal adaptation.

During probability estimation, the physical parameters remain fixed while
CE adaptation changes the scenario proposal. Appendix~\ref{subsec:surrogate_screening_integration}
records the calibrated parameters and their use in the estimator. Each
sampled scenario is evaluated with the specified dynamics and overload
criterion before likelihood weighting is applied.

\section{WECC Study}\label{sec:wecc}
The experiments use the reduced WECC 240-bus system to assess state
transfer across fault locations and operating conditions, then evaluate
probabilistic screening with the calibrated global model.

\subsection{Reference Case and Calibration Comparisons}
\label{subsec:calibration_protocol}
\label{subsec:calibration_hierarchy}
The reduced WECC 240-bus transmission system includes conventional and
renewable generation, loads, transmission corridors, and dynamic devices.
The ParaEMT implementation supplies the network and dynamic-model inputs
\cite{xiong2023paraemt,paraemt240repository}; the case description and
one-line diagram appear in Appendix~\ref{app:wecc_details}. The topology
and dynamic parameters are fixed throughout the calibration and testing
experiments. Only the specified fault, open-phase duration, observation
time, and pre-fault operating condition change.

\runin{Global calibration data and target.}
The calibration set for $\mathrm M_1$ and $\mathrm M_2$ contains 25 faulted lines, each
evaluated at 46 pre-fault operating states, giving
$25\times46=1{,}150$ training scenarios.
The operating states are generated by uniformly scaling the nominal
active-power injection vector,
\begin{equation}
    P(z_{\chi})
    =
    \chi P^{(0)},
    \qquad
    \chi\in[0.75,1.20],
    \label{eq:operating_scaling}
\end{equation}
with the steady-state equilibrium recomputed for each scale. The nominal
case is $\chi=1$; the family spans injection changes from $-25\%$ to
$+20\%$ while retaining the electrical and dynamic-model parameters
\cite{kundur1994power,chiang2011direct}.

Every calibration run, including the additional runs used for $\mathrm M_3$ and
$\mathrm M_4$, uses $\Top=\Tev=0.8\,\mathrm{s}$. The objective
in \eqref{eq:calibration_objective} therefore matches the projected
ParaEMT state at the end of the open-phase interval. This state supplies
the initial condition for subsequent restored-network evolution. The
endpoint-local mask $q=0$ selects the observations used for fitting;
the covariance target is constructed from the 46 training operating
states at the same $0.8$-s snapshot. No later testing time enters the
global calibration objective.

\runin{Line-disjoint testing of the global models.}
Let $\mathcal{E}_{\rm train}$ and $\mathcal{E}_{\rm test}$ denote the
calibration and testing fault-line sets for $\mathrm M_1$ and $\mathrm M_2$. Each contains 25 distinct
lines, and $\mathcal{E}_{\rm train}\cap\mathcal{E}_{\rm test}=\varnothing$.
The exact line identifiers and their training/test assignments are
specified in the project repository~\cite{almada_powergrid_repository}.
The split is by faulted-line identity, not a random division of the
1,150 training scenarios. Testing uses a separate $25\times46$ family
on the same prescribed operating scales. The global-model records are
therefore partitioned as
\begin{equation}
    \mathcal{D}
    =
    \mathcal{D}_{\mathrm{train}}
    \cup
    \mathcal{D}_{\mathrm{test}},
    \qquad
    \mathcal{D}_{\mathrm{train}}
    \cap
    \mathcal{D}_{\mathrm{test}}
    =
    \varnothing.
    \label{eq:data_split}
\end{equation}
For the global fits $\mathrm M_1$ and $\mathrm M_2$, only $\mathcal D_{\rm train}$ enters
parameter identification. The fitted coefficients are then fixed when
assessing the 25 previously unused faulted lines. The shared operating
scales make this a test of transfer to different fault locations within
the stated operating-state family, not a test on unseen networks or on
an independently held-out range of injection scales.

\runin{Scenario-specific calibration.}
The reference collection is enlarged for $\mathrm M_3$ and $\mathrm M_4$ by including
reclosure observations for the 25 lines in $\mathcal{E}_{\rm test}$ at
all 46 operating states, again with
$\Top=\Tev=0.8\,\mathrm{s}$. These are additional calibration runs,
not the $\Top=1.5\,\mathrm{s}$, $\Tev=2.5\,\mathrm{s}$ observations
used for assessment. Thus, a line reserved for testing $\mathrm M_1$ and $\mathrm M_2$
can supply separate fitting observations for the scenario-specific
models without entering the global fit.

Model $\mathrm M_3$ identifies 25 parameter pairs $(\alpha_e,\sigma_e)$,
one for each assessment line, using its 46 calibration operating states
jointly. Each pair is then shared across those operating states.
Model $\mathrm M_4$ identifies $25\times46=1{,}150$ pairs
$(\alpha_{e,z},\sigma_{e,z})$, one for each line and operating state.
The two specialized models use the same additional reference family
but differ in parameter sharing. Their coefficients are fitted at the
$0.8$-s target and held fixed during the subsequent assessment.

\runin{Post-reclosure assessment.}
Figure~\ref{fig:boxplots} uses $\Top=1.5\,\mathrm{s}$ and
$\Tev=2.5\,\mathrm{s}$ for the 25 testing lines across the 46 scaled
operating states. This is a deliberately demanding transfer test: the
open-phase interval is longer than in calibration, and the observation
is made one second after reclosure, while the response can still exhibit
substantial oscillations. It measures how well the calibrated description
reproduces the subsequent electromechanical state, rather than how
closely it fits the training snapshot. These settings are not asserted
to be a worst case over all faults or operating conditions. The full
transient in Fig.~\ref{Point3}(a) illustrates oscillatory post-reclosure
behavior under its own, shorter $0.5$-s open-phase interval; it is not
the trajectory used to define the timing of Fig.~\ref{fig:boxplots}.

Within every reference--surrogate pair, the initial operating state,
faulted line, switching times, phase reference, and observation operator
are matched. All five models are scored on the same network-wide
observations, with no refitting to the assessment response.
For $\mathrm M_1$ and $\mathrm M_2$, the assessment combines transfer to unseen lines
with transfer in duration and observation time; for $\mathrm M_3$ and $\mathrm M_4$,
it measures the latter transfer after case-specific calibration.
Subsequent probability estimates use only $\mathrm M_2$, with its global pair
held fixed.

\runin{Separate temporal checks.}
Appendices~\ref{subsec:fault_duration} and~\ref{subsec:dead_time} vary
$\Tev$ and $\Top$, respectively, with unscaled injections and the global
parameter pair held fixed. These experiments separate observation-time
sensitivity from switching-duration sensitivity.

\begin{figure*}[t]
\centering
\includegraphics[width=\textwidth]{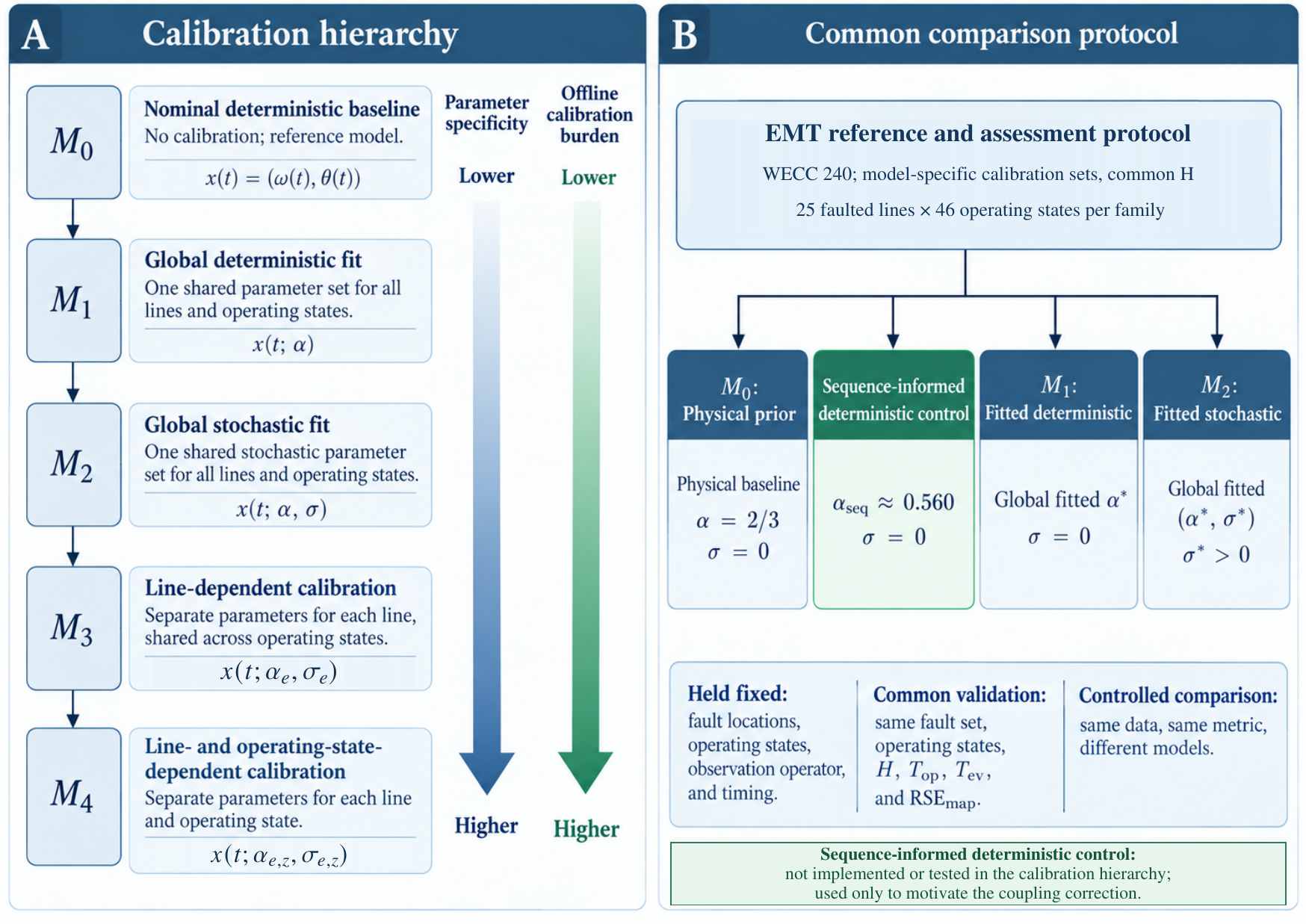}
\caption{\textbf{Calibration hierarchy and common comparison protocol.}
\textbf{(A)} $\mathrm M_0$ fixes $(\alpha,\sigma)=(2/3,0)$; $\mathrm M_1$ fits one
global deterministic coupling and $\mathrm M_2$ one global stochastic pair.
$\mathrm M_3$ fits 25 line-specific pairs and $\mathrm M_4$ fits 1,150
line--operating-state-specific pairs. All fits use
$\Top=\Tev=0.8\,\mathrm{s}$. The specialized models additionally use
reference data from the 25 assessment lines, which are excluded from the
global fits. The arrows describe fitting specificity and effort, not
measured runtime ratios.
\textbf{(B)} The common protocol fixes the WECC model, observations, and
assessment conditions. Figure~\ref{fig:boxplots} scores every model on
25 lines and 46 operating states at $(\Top,\Tev)=(1.5,2.5)\,\mathrm{s}$.
Those line identities are unseen by $\mathrm M_1/\mathrm M_2$, but have separate
$0.8$-s calibration records for $\mathrm M_3/\mathrm M_4$. The $25\times46$ count
identifies a reference family, not one pool split among all models.
The sequence-informed deterministic reference only motivates the coupling
correction; it is not implemented or tested in the hierarchy. The line
and operating-state indices are $e$ and $z$, respectively.}
\label{fig:sequence_control_design}
\end{figure*}

\runin{Comparison hierarchy.}
Figure~\ref{fig:sequence_control_design}(A) defines the five model levels.
The baseline $\mathrm M_0$ fixes $\alpha=2/3$ and $\sigma=0$; $\mathrm M_1$ fits a
shared deterministic coupling, and $\mathrm M_2$ jointly fits the shared coupling
and stochastic amplitude. Models $\mathrm M_3$ and $\mathrm M_4$ add line and operating-state
dependence. The comparison tests increasing parameter specificity within this model family.

Panel (B) also shows the sequence-informed deterministic reference
$\alpha_{\rm seq}\simeq0.560$, $\sigma=0$. This reference motivates
the coupling correction; it is neither implemented nor tested as part
of the calibration hierarchy. The
common protocol means matching all model comparisons within each stage;
it does not mean using the same fault catalogue or time settings for
training and testing.

\subsection{Probabilistic Screening Setup}
\label{subsec:wecc_screening_setup}
CEM and MC use the same calibrated surrogate, nominal contingency law,
screening horizon $\Thor=10\,\mathrm{s}$, and trajectory-evaluation
interval $\Delta t=0.01\,\mathrm{s}$ on the reduced WECC system
\cite{price2011reduced}. The tested global thresholds are
$\gamma\in\{0,0.5,5,10\}\,\mathrm{s}$ in
\eqref{eq:global_probability}. The global parameters remain fixed during
sampling as specified in \eqref{eq:screening_parameters}.

The nominal law selects among $n_{\rm E}=329$ modeled faulted lines uniformly and
uses an exponential open-phase duration with mean
$\tau_0=0.1\,\mathrm{s}$; see \eqref{eq:cem_nominal}.
The probabilities are conditional on a fault occurring within this line
set. This 329-line screening catalogue is broader than either of the
25-line calibration and testing sets. Applying the global pair throughout that catalogue extrapolates it beyond
the EMT-assessed line sets. The nominal mean $0.1\,\mathrm{s}$ specifies the probability
law; it is distinct from the fixed $0.8$-s calibration duration and
$1.5$-s hierarchy-test duration. CE adaptation uses $250$ trajectories per iteration for $20$
iterations, with the $25$ largest scores forming the elite set. These
$5{,}000$ adaptation evaluations are separate from the final estimation
batch and are included in the reported total CEM runtimes.

The duration-resolved assessment in Figs.~\ref{fig2} and \ref{fig3}
uses $N=150{,}000$ surrogate simulations. It varies the open-phase duration
rather than marginalizing it under the nominal duration law. Its
line-specific probabilities must therefore be distinguished from the
duration-marginalized map in Fig.~\ref{fig:contingency_risk_heatmap}.
The event definitions themselves are those of
Sec.~\ref{subsec:dynamic_overload}; they are not redefined for the WECC
results.

\section{Results}\label{sec:results}
\label{sec:calibration_validation}
\label{sec:screening}\label{sec:deployment}
The results address three questions: what the surrogate retains during
the open-phase interval, how much parameter flexibility is needed for
state transfer, and what is learned by resolving overload risk across
initiating faults and monitored elements. Sampling diagnostics are
reported separately from model accuracy.

\subsection{Open-Phase Trajectory Comparison}
\label{subsec:open_phase_comparison}
\begin{figure*}[t]
\centering
\includegraphics[width=\textwidth]{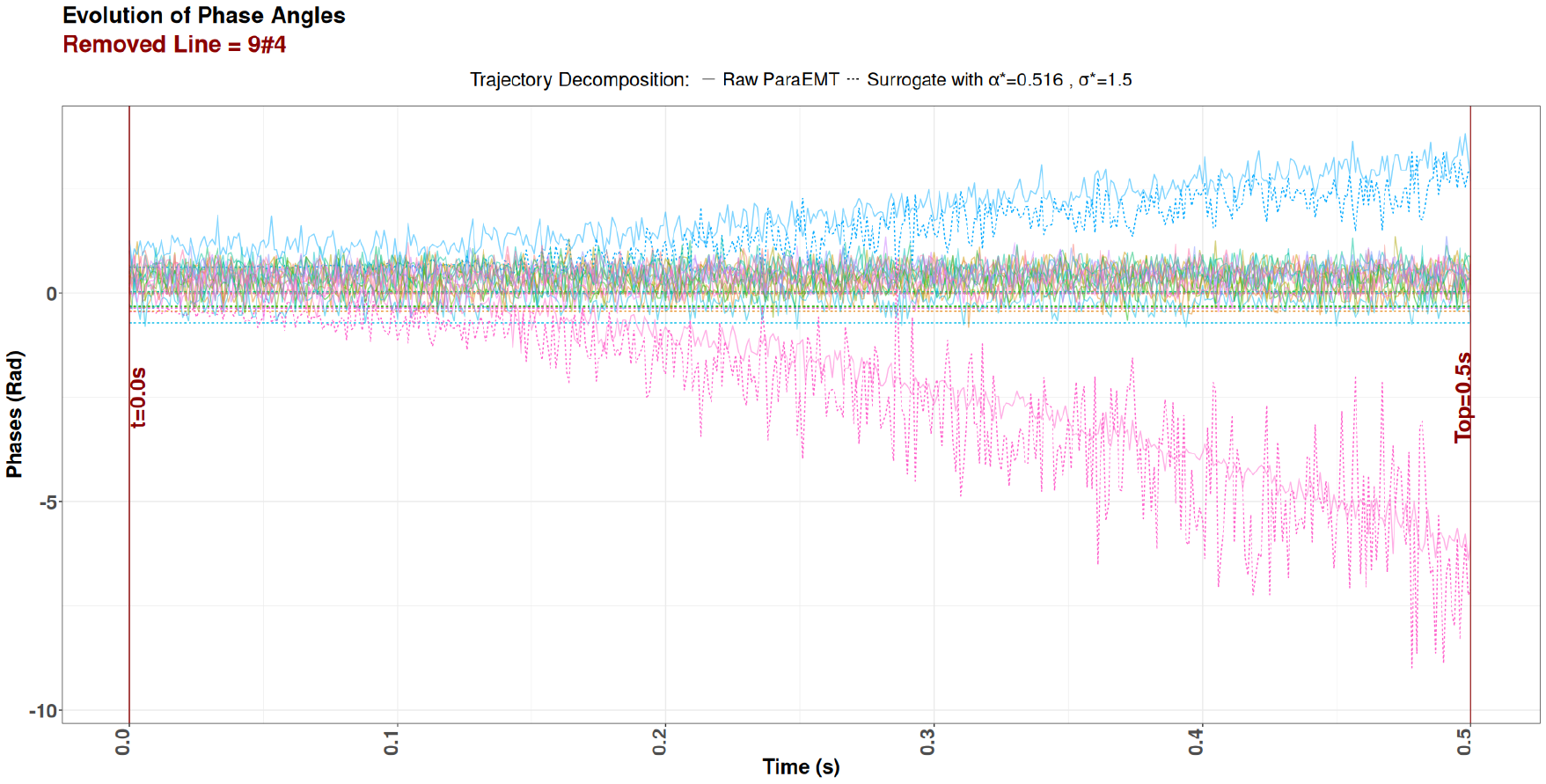}
\caption{\textbf{ParaEMT--surrogate comparison during the open-phase interval.}
Phase-angle trajectories for the illustrative large-excursion $9\#4$
contingency over $0\leq t\leq\Top=0.5\,\mathrm{s}$.
Solid curves show ParaEMT; dashed curves show one realization of
\eqref{eq:surrogate_sde} with $(\alpha^*,\sigma^*)=(0.516,1.5)$,
not an ensemble mean or confidence band. The vertical markers identify
pole opening and reclosure. Figure~\ref{Point3} continues this same
example to $10\,\mathrm{s}$ and displays the corridor separation
$|\theta_9-\theta_4|$. The full-trajectory endpoint is not the
$0.8$-s calibration target or the $2.5$-s hierarchy-test time.}
\label{fig:open_phase_emt_surrogate}
\end{figure*}

Figure~\ref{fig:open_phase_emt_surrogate} compares ParaEMT and
one realization of the calibrated stochastic surrogate for the
illustrative large-excursion $9\#4$ contingency during $[0,\Top]$,
with $\Top=0.5\,\mathrm{s}$. Both calculations show an uneven phase
response: many traces remain in a relatively narrow range, while others
develop larger excursions before reclosure. The surrogate follows this
overall separation of responses, with differences in fluctuation
amplitudes. The comparison
makes the retained response and the remaining trajectory discrepancies
visible during the interval in which the modified coupling and stochastic
forcing act.

The full trajectory and corridor-separation diagnostic in
Fig.~\ref{Point3} continue this example through reclosure. The hierarchy
below provides the quantitative comparison of projected mean states.

\subsection{Calibration Accuracy and Parameter Dependence}
\label{subsec:global_wecc_calibration}
\label{subsec:calibration_interpretation}
Joint optimization of the global stochastic model $\mathrm M_2$ under the
protocol in Sec.~\ref{subsec:calibration_protocol} gives
\begin{equation}
    (\alpha^{\ast},\sigma^{\ast})
    =
    (0.516,\ 1.5).
    \label{eq:global_optimum}
\end{equation}
At $\Top=1.5\,\mathrm{s}$ and $\Tev=2.5\,\mathrm{s}$, the fixed global
pair gives $\mathrm M_2$ errors of approximately $0.0968$--$0.15$ across the
25 testing lines and 46 uniformly scaled operating states, with mean
$0.098$. The one-sided distribution places the mean near its lower end.

\begin{figure*}[t]
\centering
\makebox[\textwidth][c]{\includegraphics[width=18cm, height=10cm]{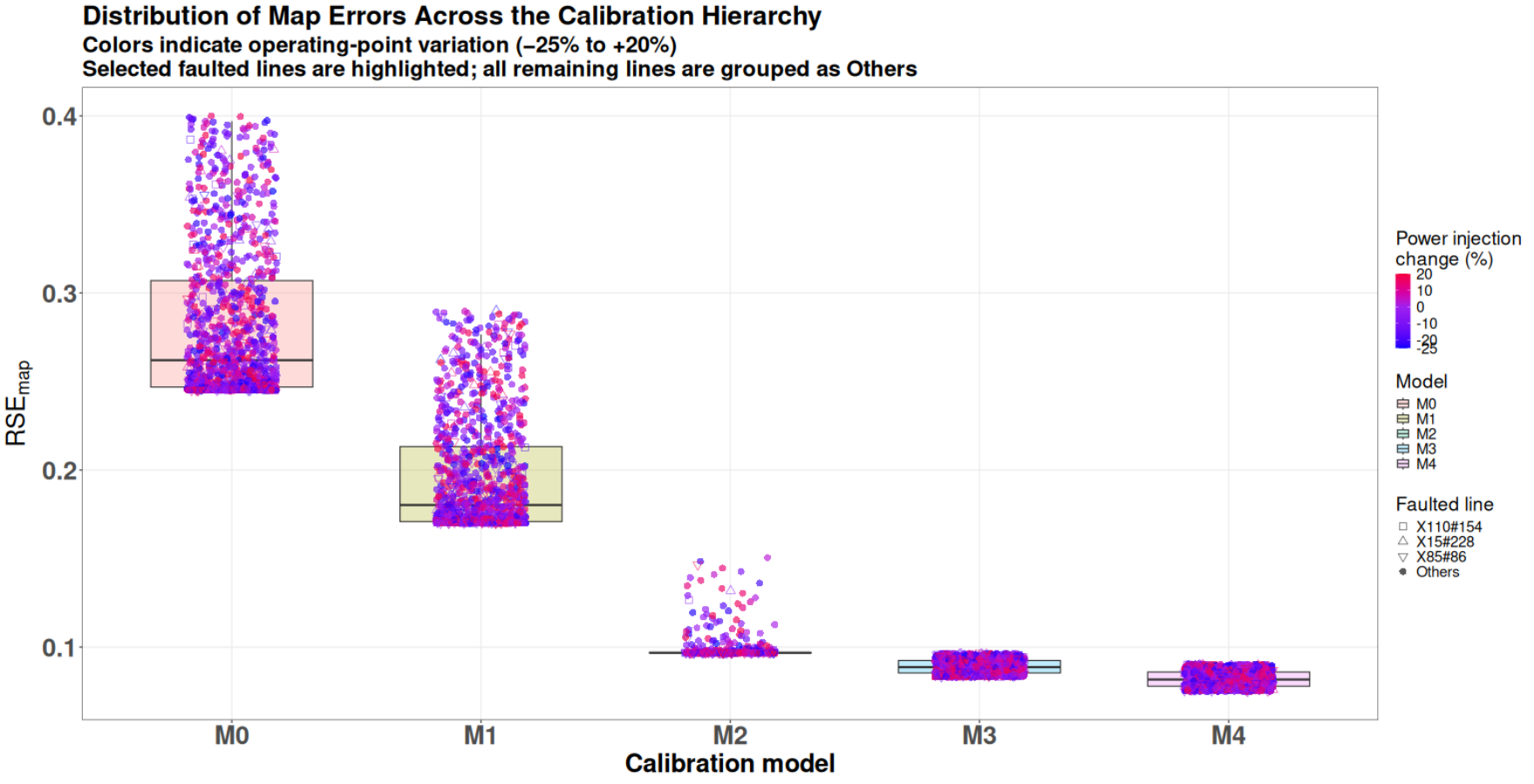}}
\caption{Distribution of $RSE_{\rm map}$ across the WECC 240-bus model
hierarchy at $\Top=1.5\,\mathrm{s}$ and $\Tev=2.5\,\mathrm{s}$.
All models are scored on the same 25 lines and 46 uniformly scaled
operating states. Marker shapes identify lines; colors identify injection
changes. $\mathrm M_1/\mathrm M_2$ were fitted on 25 disjoint calibration lines.
$\mathrm M_3$ instead uses 25 line-specific pairs fitted across the 46 operating
states, and $\mathrm M_4$ uses 1,150 line--state-specific pairs, obtained from
additional reference runs on the assessment lines. Every fit uses
$\Top=\Tev=0.8\,\mathrm{s}$ and is held fixed here. Thus, the
$\mathrm M_3/\mathrm M_4$ results measure duration and observation-time transfer after
case-specific calibration, not transfer to unseen lines.
Their smaller additional gains beyond $\mathrm M_2$ require more separate fits.
Unlike Figs.~\ref{fig:fault_duration}--\ref{fig:dead_time}, this
comparison includes injection scaling.}
\label{fig:boxplots}
\end{figure*}

Figure~\ref{fig:boxplots} compares the hierarchy under this common
assessment. The fixed deterministic baseline $\mathrm M_0$ has errors of
approximately $0.24$--$0.4$, while fitting the deterministic coupling
in $\mathrm M_1$ brings most cases into the $0.17$--$0.29$ range. The dominant
reduction occurs from $\mathrm M_1$ to $\mathrm M_2$. Because
$\mathrm{RSE}_{\rm map}$ scores the projected mean, this comparison
shows the contribution of the noise-induced drift to state transfer
within the tested hierarchy.

Models $\mathrm M_3$ and $\mathrm M_4$ produce narrower error distributions and smaller
further reductions using the additional case-specific calibration records
specified in Sec.~\ref{subsec:calibration_protocol}. Substantial variation
in their fitted parameters, reported in \eqref{eq:parameter_ranges} and
Fig.~\ref{fig:opsensitivity}, therefore accompanies relatively modest
accuracy gains. This is the main calibration trade-off: the global pair
captures most of the improvement, while the specialized fits require
additional reference observations and more independent parameter fits.
The subsequent screening calculations retain $\mathrm M_2$ on this basis.

\runin{Coupling and temporal transfer.}
The fitted coupling lies below the phase-counting baseline, in the same
direction as the representative sequence-network estimate
$\alpha_{\rm seq}\simeq0.560$. The residual difference distinguishes the
line-level sequence calculation from the joint fit to the finite-time
network response.

The unscaled-injection studies in
Appendices~\ref{subsec:fault_duration} and~\ref{subsec:dead_time}
separate two temporal effects. After a fixed reclosure, discrepancies
continue to increase even though the stochastic forcing has ceased:
restoring the nominal operator does not erase the mismatch in the attained
state. Varying the open-phase duration instead changes how long the
modified dynamics act before that state is transferred. These trends
complement the hierarchy comparison by showing how error depends on time
as well as on the fitted parameterization.

\subsection{Overload Probabilities and Sampling Efficiency}
\label{subsec:wecc_probabilistic}
We evaluate the probabilities defined in
\eqref{eq:global_probability} using the WECC setup in
Sec.~\ref{subsec:wecc_screening_setup}.

Table~\ref{tab:reduced_data_results} reports the probability estimates,
final estimation-batch sizes, and total runtimes. Including the adaptation
budget specified in Sec.~\ref{subsec:wecc_screening_setup}, CEM uses
$30{,}000$--$160{,}000$ surrogate evaluations, compared with
$125{,}000$--$300{,}000$ for MC. Tables~\ref{tab:wecc_sampling_diagnostics}
and~\ref{table4} in Appendix~\ref{app:wecc_sampling_details} provide the
sampling diagnostics and stage-by-stage costs.

\begin{table}[t]
\caption{Probability estimates and total runtimes for CEM and Monte Carlo
(MC) \cite{kroese2012monte}. Here, $N_{\rm final}$ is the final
estimation-batch size. CEM total runtimes include both the $5{,}000$
adaptation evaluations and final probability estimation.}
\label{tab:reduced_data_results}
\small
\centering
\begin{tabular}{llccc}
\hline
$\boldsymbol{\gamma}$ \textbf{[s]} & \textbf{Method} & $N_{\rm final}$ &
$\boldsymbol{\widehat Q}$ & \textbf{Total [s]} \\
\hline
$\boldsymbol{0}$
& \textbf{CEM}
& $2.500 \times 10^{4}$
& $\boldsymbol{5.600 \times 10^{-2}}$
& $\boldsymbol{4.462 \times 10^{3}}$ \\
& \textbf{MC}
& $1.250 \times 10^{5}$
& $\boldsymbol{5.606 \times 10^{-2}}$
& $\boldsymbol{2.210 \times 10^{4}}$ \\
\hline
$\boldsymbol{0.5}$
& \textbf{CEM}
& $4.000 \times 10^{4}$
& $\boldsymbol{2.628 \times 10^{-2}}$
& $\boldsymbol{7.177 \times 10^{3}}$ \\
& \textbf{MC}
& $1.750 \times 10^{5}$
& $\boldsymbol{2.636 \times 10^{-2}}$
& $\boldsymbol{3.185 \times 10^{4}}$ \\
\hline
$\boldsymbol{5.0}$
& \textbf{CEM}
& $1.050 \times 10^{5}$
& $\boldsymbol{7.743 \times 10^{-3}}$
& $\boldsymbol{1.970 \times 10^{4}}$ \\
& \textbf{MC}
& $2.500 \times 10^{5}$
& $\boldsymbol{7.872 \times 10^{-3}}$
& $\boldsymbol{4.536 \times 10^{4}}$ \\
\hline
$\boldsymbol{10.0}$
& \textbf{CEM}
& $1.550 \times 10^{5}$
& $\boldsymbol{2.723 \times 10^{-3}}$
& $\boldsymbol{2.828 \times 10^{4}}$ \\
& \textbf{MC}
& $3.000 \times 10^{5}$
& $\boldsymbol{2.900 \times 10^{-3}}$
& $\boldsymbol{5.547 \times 10^{4}}$ \\
\hline
\end{tabular}
\end{table}

Both estimators resolve the decrease in exceedance probability from
approximately $0.056$ at $\gamma=0$ to $0.0077$--$0.0079$ at
$\gamma=5\,\mathrm{s}$ and $0.0027$--$0.0029$ at
$\gamma=10\,\mathrm{s}$. The absolute differences between their point
estimates range from $6.0\times10^{-5}$ to $1.77\times10^{-4}$.
The estimates distinguish any positive overload from larger accumulated
network exposure and reach the $10^{-3}$ probability range.

At every reported threshold, CEM uses a smaller final batch and less
total runtime than MC. The MC-to-CEM runtime ratios are $4.952$,
$4.438$, $2.303$, and $1.961$ for increasing $\gamma$. The advantage
at these budgets is therefore largest for the lowest threshold and
narrows as the target event becomes rarer. At
$\gamma=10\,\mathrm{s}$, CEM reaches an estimate of $0.002723$ in
$28{,}281.6\,\mathrm{s}$ including adaptation, compared with
$0.002900$ in $55{,}468.7\,\mathrm{s}$ for MC.

The contribution-variance analysis in
Appendix~\ref{app:wecc_sampling_details} favors CEM at
$\gamma=0.5$ and $5\,\mathrm{s}$, and MC at $0$ and
$10\,\mathrm{s}$. Adaptation selects the largest global scores, so its per-sample benefit
varies with the final exceedance threshold.

\subsection{Contingency-to-Vulnerable-Element Map}
\label{subsec:risk_ranking}
The conditional probabilities defined in
Sec.~\ref{subsec:dynamic_overload} resolve the location of the consequences
associated with each initiating fault. They are conditional on that fault
and the specified operating state, rather than estimates of the fault's
occurrence frequency.

\begin{figure*}[t]
\centering
\includegraphics[width=\textwidth]{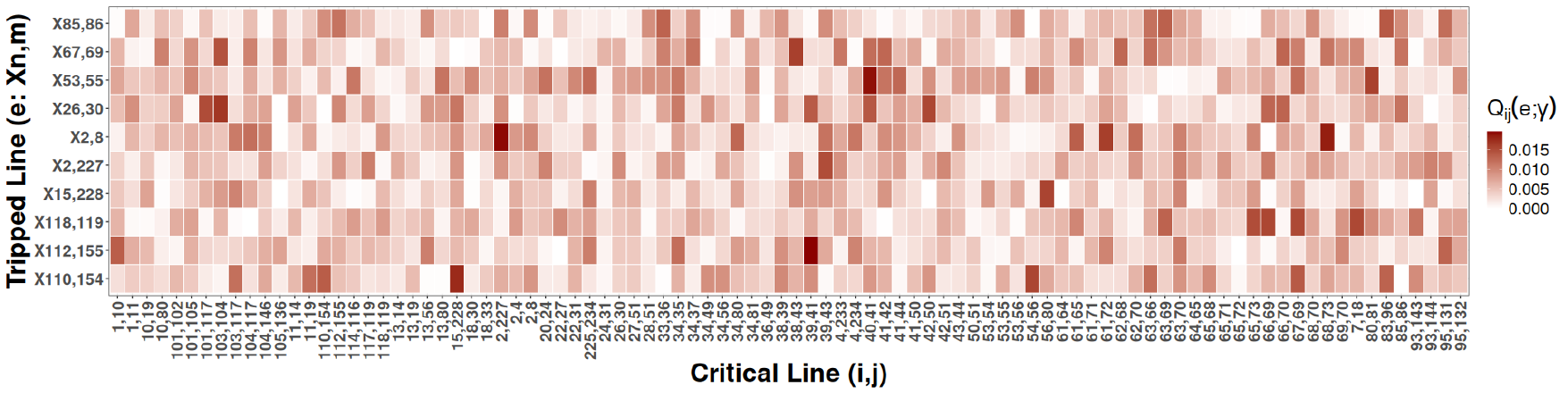}
\caption{Contingency-to-vulnerable-element map for the WECC 240-bus
assessment at $\gamma=1\,\mathrm{s}$ over $\Thor=10\,\mathrm{s}$.
The x-axis shows critical transformer-connected monitored lines; the
y-axis identifies the ten displayed initiating faults. The comma-separated
axis label $Xn,m$ denotes the same fault identifier as \texttt{Xn\#m}
in the experiment records. The row order is not a scalar-risk ranking. Each cell represents
$Q_{ij}(e;\gamma)=P\left[S_{ij}(\theta^{e})>\gamma\right]$, the probability of more than one second of accumulated
overload on element $\{i,j\}$ conditional on fault $e$, with open-phase
duration and surrogate uncertainty marginalized under the specified
nominal law. Entries with $\{i,j\}\neq e$ identify secondary exposure.
The continuous colorbar includes labeled values through $0.015$;
its colors are not the duration-resolved warning classes.}
\label{fig:contingency_risk_heatmap}
\end{figure*}

Figure~\ref{fig:contingency_risk_heatmap} presents the conditional map
at $\gamma=1\,\mathrm{s}$ for ten initiating faults and the monitored
elements shown. Rows are initiating faults, columns are monitored
elements, and each cell is $Q_{ij}(e;1\,\mathrm{s})$. The row order
is not a descending ranking of the scalar scores $Q(e;\gamma)$.
The figure instead resolves the fault--element pairs underlying the
screening output.

The colorbar has labeled ticks at $0$, $0.005$, $0.010$, and $0.015$.
Larger entries occur in different columns for different faults.
The \textbf{X2\#227} contingency, for example, has a prominent
secondary-risk entry on $\{61,65\}$. A prominent entry for
\textbf{X112\#155} instead lies on the faulted line $\{112,155\}$
itself. These contrasting cases show why a global score alone is
insufficient: it does not identify the monitored element or distinguish
self-exposure from secondary consequences.

\runin{Interpretation of the screening thresholds.}
The duration-resolved assessment in Fig.~\ref{fig3} uses
$\gamma=1\,\mathrm{s}$ and probability boundaries
$0.025$ and $0.04$: \textbf{safe} denotes probabilities below
$0.025$, \textbf{warning} values from $0.025$ through $0.04$, and
\textbf{emergency} values above $0.04$.
The heatmap marginalizes open-phase duration under the nominal law,
whereas the curves in Appendix~\ref{app:screening_details} condition
on duration. Its continuous color scale is independent of those classes.
Section~\ref{subsec:scope_limitations} states the interpretation of
these screening thresholds.

The duration-resolved curves supply complementary information: an element
can have substantial exposure under a long open-phase interval while
retaining a small duration-marginalized probability under the nominal law.

\subsection{Computational Cost and Reuse of Calibration}
\label{subsec:surrogate_vs_emt}

Replacing repeated EMT simulations with reduced trajectories addresses
the largest per-realization computational burden. A direct EMT ensemble
requires a high-fidelity trajectory for each sampled contingency and
open-phase duration. The reported mean surrogate evaluation time is
approximately $0.19\,\mathrm{s}$; the published WECC benchmark reports
$29\,\mathrm{s}$ on one core and $28\,\mathrm{s}$ on eight cores
for a $1$-s ParaEMT simulation \cite{xiong2023paraemt}. The benchmark
settings are given in Table~\ref{tab:paraemt_published_benchmark}.
Section~\ref{subsec:scope_limitations} states the scope of this
computational comparison.

Reference generation and parameter identification are offline costs.
The screening ensemble reuses the calibrated $\mathrm M_2$ model for every
adaptation and estimation trajectory. This reuse is shared by CEM and MC;
their sampling costs and precision are compared in
Sec.~\ref{subsec:wecc_probabilistic}. The additional calibration records
used by $\mathrm M_3/\mathrm M_4$ belong to the hierarchy assessment, not to the repeated
screening calculation.

\section{Discussion}\label{sec:conclusion}\label{sec:discussion}
The advance over the earlier screening construction is the identification
of its effective parameters from EMT and the assessment of their
transferred response. The WECC hierarchy supports a compact global
parameterization, while the specialized fits show the value of additional
case-specific information. The temporal studies further identify how
switching duration and subsequent propagation affect the state mismatch.

Deterministic moment propagation separates parameter identification from
trajectory sampling. The calibrated model can then be reused across
counterfactual scenarios, with CE adaptation directing evaluations toward
large accumulated exposure. The resulting fault--element map identifies
where detailed high-fidelity follow-up is most valuable, connecting EMT
reference information to probabilistic screening.

\subsection{Scope and limitations}
\label{subsec:scope_limitations}
The assessment covers the reduced WECC network, the prescribed uniform
injection-scale family, and the timing protocols in
Sec.~\ref{subsec:calibration_protocol}. The fitted coefficients are
effective parameters for those conditions. The stochastic correction is
an empirical representation of unresolved variability, rather than a
white-noise model derived from the EMT equations. The hierarchy establishes
its value within the tested family; distinguishing this mechanism from
other reduced descriptions with comparable fitting flexibility requires
additional comparisons. Applying the global pair to the broader
329-line screening catalogue extends beyond the EMT-assessed line sets.
Changes in topology, device models, or operating conditions require
further assessment of the calibration domain.

The metric $\mathrm{RSE}_{\rm map}$ evaluates the projected mean state
at a specified time. It does not establish agreement of the complete
waveform, the path-dependent overload duration, or the associated tail
probability. Errors that change whether $S>\gamma$ can affect a rare-event
probability substantially, whereas state discrepancies away from the
threshold may have little effect on that event. CEM and MC evaluate the
same surrogate: their agreement assesses probability estimation under
that model, not independent agreement with EMT. Increasing the sampling
budget reduces statistical uncertainty without correcting model
mismatch. EMT-informed overload, probability, and ranking comparisons
under the same scenario law remain the corresponding physical-validation
tasks. The state-error ordering of the hierarchy is therefore not an
established ordering of probability accuracy. The second-moment spectral
diagnostic likewise concerns the reduced stochastic dynamics.

The overload indicator uses a linearized power--angle relation and
measures accumulated threshold exposure. The large corridor separations
in Fig.~\ref{Point3} limit its small-angle interpretation during that
transient. Outside this regime, the indicator serves as a reduced-model
severity measure for high-fidelity follow-up, rather than an exact AC
power-flow or thermal-damage calculation. The duration threshold,
probability-class boundaries, and angle-diagnostic band are illustrative;
they are not NERC or WECC criteria or relay settings. In particular,
\textbf{safe} denotes the low-probability screening class, not a guarantee
of security. The resulting map prioritizes assessment rather than
prescribing protection actions.

The reported probabilities also depend on the specified exponential
open-phase-duration law with mean $0.1\,\mathrm{s}$ and on successful
reclosure. This law is a duration model for the study, not an empirical
distribution of protection outcomes. A failed reclosure requires a
different switching sequence. Protection-specific branches and broader
operating-state families are extensions of the present formulation.

The runtime comparison uses the reported sample budgets and their
associated statistical uncertainties, rather than a matched-precision
benchmark. Its contribution-variance advantage is event-dependent.
Calibration fit counts describe fitting requirements, not measured
wall-clock multipliers. The published ParaEMT timings provide external
computational context: the model fidelity, numerical settings, and
potential hardware differences preclude a controlled speedup claim or an
EMT-ensemble runtime extrapolation from those values. The present
screening implementation takes hours per threshold and remains a proof
of concept. Complete update latency must include adaptation and final
estimation at the required precision. Parallel trajectory evaluation and
reuse of a CE proposal across nearby operating states remain options to
test before operational deployment.

\section*{Data Availability}
The project code and WECC calibration materials are available through
the PowerGridREserach repository~\cite{almada_powergrid_repository}
and its calibration-data directory~\cite{almada_wecc_calibration_data}.
The model-specific calibration and assessment allocations are specified
in Sec.~\ref{subsec:calibration_protocol}.
Historical Israel-case materials are stored separately. ParaEMT is
available from its upstream repository~\cite{paraemt_public_repository}.

\bibliographystyle{ieeetr}
\bibliography{refs}

\appendix

\section{EMT Reference, WECC Case, and Event Notation}\label{app:reference_details}

This appendix identifies the reference implementation and case inputs
used in Sec.~\ref{sec:wecc} and summarizes the switching and observation
settings using the time conventions in Table~\ref{tab:global_notation}.

\subsection{Reference Model and WECC Case}
\label{app:wecc_details}\label{subsec:paraemt}
ParaEMT couples three-phase nodal network equations with synchronous-machine,
excitation, governing, stabilizer, and converter models. Its
differential-algebraic formulation is evaluated using trapezoidal
integration and sparse linear algebra
\cite{xiong2024open,ma2023basics,talpur2017testing}. The observation
operator in \eqref{eq:emt_observation} selects the electromechanical
variables compared with the surrogate; fast electrical and device states
remain in the reference model.

\begin{figure}[t]
    \centering
    \includegraphics[width=\columnwidth]{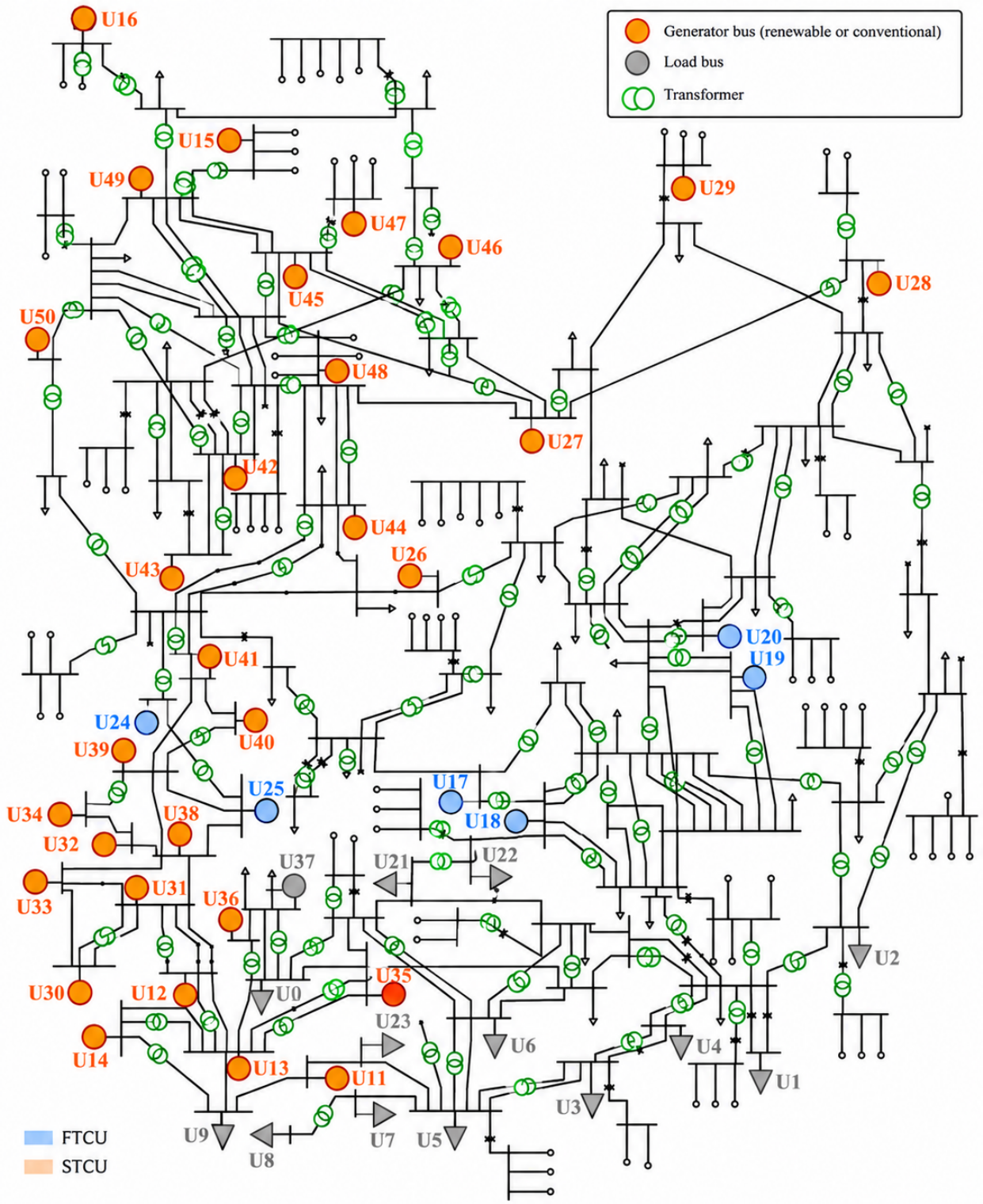}
    \caption{One-line representation of the reduced WECC 240-bus calibration
    case, adapted from \cite{wang2025comprehensive}.}
    \label{fig:wecc240}
\end{figure}

The experiments use the public reduced WECC 240-bus case, with conventional
and renewable generation, loads, transmission corridors, and approximately
20\% inverter-based-resource penetration \cite{xiong2023paraemt}.
The case is distributed in the \texttt{models/240bus\_psse} directory
of ParaEMT, with the \texttt{WECC240.raw} network description and
\texttt{240bus.xls} model data \cite{paraemt240repository}.
All calibration levels use this same topology and set of dynamic-model
parameters; only the specified contingency and operating condition change.
Figure~\ref{fig:wecc240} shows the network.

\subsection{Timeline Reference}\label{app:timeline_reference}
The switching sequence is defined in Sec.~\ref{subsec:terminology}
and illustrated in Fig.~\ref{fig:overview}. Table~\ref{tab:global_notation}
collects the time conventions; the experiment-specific settings and
calibration allocations are given in
Sec.~\ref{subsec:calibration_protocol}. The $9\#4$ illustration uses
$\Top=0.5\,\mathrm{s}$ and ends at $10\,\mathrm{s}$. These settings
remain distinct from the nominal screening-duration distribution.

\section{Sequence-Network Interpretation of the Coupling}\label{app:sequence_network}

A quantitative sequence-network calculation provides a first physical check
of the fitted coupling. Among the WECC testing contingencies, line
$X110\#154$ has positive-sequence reactance $X_1=0.02886$ p.u., equal to the
median reactance of the tested line set, and is therefore selected as a
representative case. For a one-conductor-open condition, the negative- and
zero-sequence networks contribute through
\[
    Z_{20}
    =
    Z_2\parallel Z_0
    =
    \frac{Z_2Z_0}{Z_2+Z_0},
\]
so that the effective positive-sequence transfer impedance is
\[
    Z_{\rm op}
    =
    Z_1+(Z_2\parallel Z_0).
\]
Using the standard transmission-line approximation $Z_2\simeq Z_1$ together
with representative overhead-line sequence ratios
$R_0/R_1=13.03$ and $X_0/X_1=3.49$ gives
\[
    Z_1
    \simeq
    0.00202+j0.02886,
    \qquad
    Z_0
    \simeq
    0.02628+j0.10072,
\]
and
\[
    Z_2\parallel Z_0
    \simeq
    0.00249+j0.02257.
\]
Because $\alpha$ scales the line coupling in
\eqref{eq:effective_coupling}, the corresponding comparison is made at the
susceptance level,
\[
    \alpha_{\rm seq}
    =
    \frac{
    \left|\operatorname{Im}
    \left\{1/[Z_1+(Z_2\parallel Z_0)]\right\}\right|
    }{
    \left|\operatorname{Im}(1/Z_1)\right|
    }
    \simeq
    0.560.
\]
The sequence-network prediction is therefore substantially below the naive
phase-counting value $\frac{2}{3}\simeq0.667$ and moves in the direction of the
EMT-calibrated value $\alpha^*=0.516$. Relative to the $\frac{2}{3}$ baseline, this
sequence correction accounts for approximately $71\%$ of the observed
reduction toward the fitted value. The residual difference,
$\alpha_{\rm seq}-\alpha^*\simeq0.044$, shows that the line-level sequence
equivalent alone does not determine the calibrated coupling.

\begin{figure}[t]
\centering
\includegraphics[width=\columnwidth]{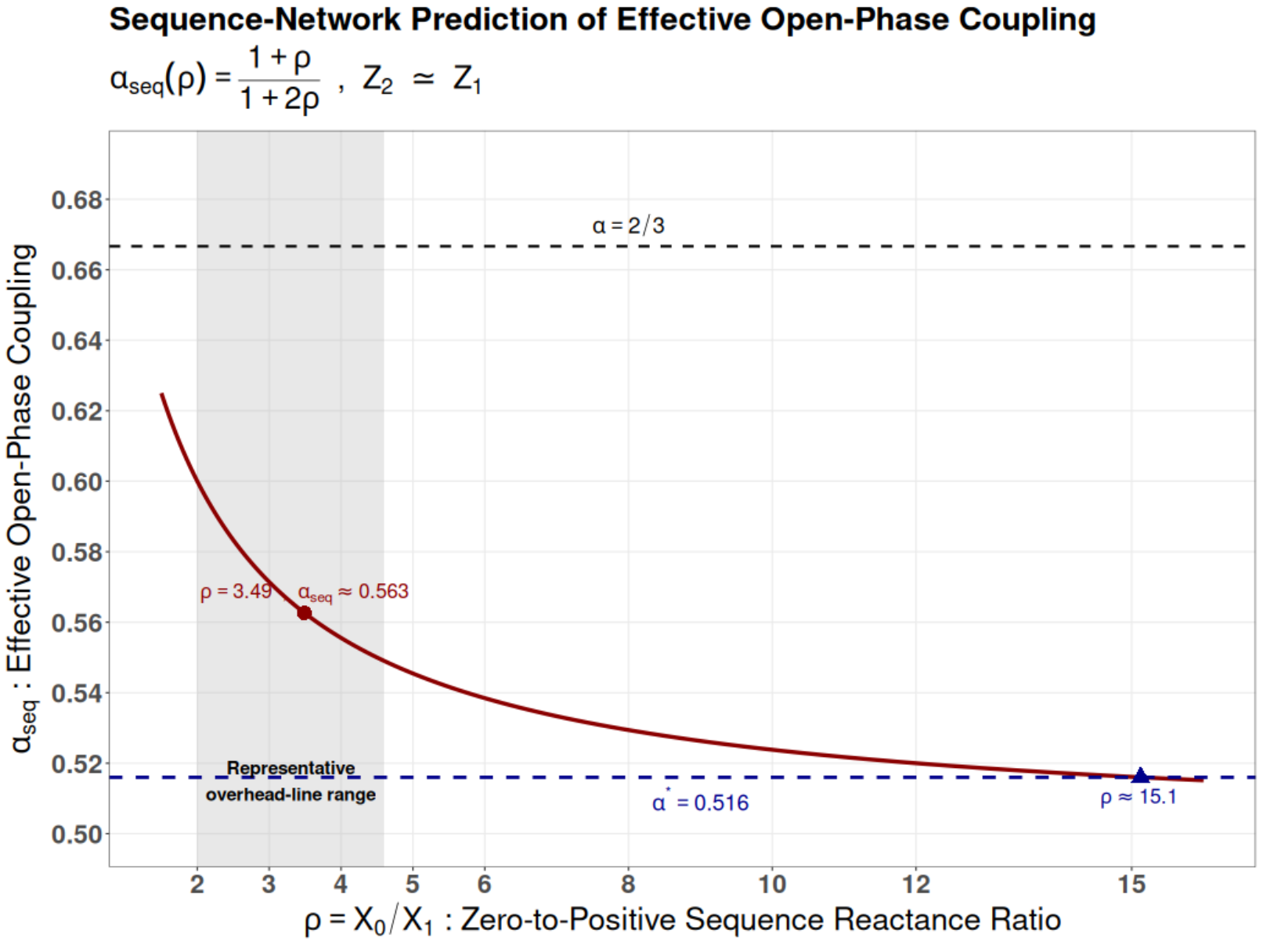}
\caption{Sequence-network coupling versus $\rho=X_0/X_1$, with
$Z_2\simeq Z_1$. Shading marks the representative overhead-line range
$2.0$--$4.6$. The representative 230-kV ratio $\rho=3.49$ gives
$\alpha_{\rm seq}\approx0.563$
(complex-impedance estimate: $0.560$); reproducing the fitted
$\alpha^*=0.516$ by this correction alone requires $\rho\approx15.1$.}
\label{fig:sequence_alpha}
\end{figure}

Figure~\ref{fig:sequence_alpha} examines the sensitivity of this result to the unavailable line-specific WECC zero-sequence parameters. Writing
$\rho=X_0/X_1$, neglecting resistance for this sensitivity calculation, and
using $Z_2\simeq Z_1$ gives
\[
    \alpha_{\rm seq}(\rho)
    =
    \frac{1+\rho}{1+2\rho}.
\]
Representative overhead-line values $\rho\simeq2.0$--$4.6$ give
$\alpha_{\rm seq}\simeq0.600$--$0.549$. The representative value
$\rho=3.49$ gives $\alpha_{\rm seq}\simeq0.563$, consistent with the
complex-impedance estimate $0.560$ above. In contrast, reproducing
$\alpha^*=0.516$ through the line-level sequence correction alone would
require $\rho\simeq15.1$, well outside the representative range considered
here.

\section{Local Calibration and Moment Diagnostics}
\label{app:calibration_details}
\subsection{Graph-Local Calibration}\label{subsec:graph_local}

Calibration uses the electromechanical response at the endpoints of the
faulted line, where the stochastic forcing in \eqref{eq:surrogate_sde}
acts \cite{pavliotis2014stochastic}. The endpoint restriction applies to
parameter identification, not to propagation through the network.

Let $\mathcal{G}=(\mathcal{V},\mathcal{E})$ denote the structure-preserving
network graph introduced in Sec.~\ref{subsec:postclearing}. For a bus
$a\in\mathcal{V}$ and a faulted line $e=(i,j)$, define the graph distance
from the bus to the line by
\begin{equation}
    d(a,e)
    =
    \min\{d(a,i),d(a,j)\},
    \label{eq:line_distance}
\end{equation}
where $d(\cdot,\cdot)$ denotes the shortest-path distance on
$\mathcal{G}$. For a locality radius $q$, define
\begin{equation}
    \mathcal{N}_{q}(e)
    =
    \left\{
        a\in\mathcal{V}:d(a,e)\leq q
    \right\}.
    \label{eq:local_neighborhood}
\end{equation}

The present analysis uses the endpoint-local case $q=0$. The general
neighborhood in \eqref{eq:local_neighborhood} then reduces to the endpoint
set introduced in \eqref{eq:fault_endpoints}:
\begin{equation}
    \mathcal{N}_{0}(e)=\{i,j\}.
\end{equation}
Only these two buses enter the calibration objective. Broader neighborhoods,
including $q=1$, are not used in the reported fits.

Let $W_e^{(q)}$ denote the corresponding diagonal bus mask,
\begin{equation}
    W_e^{(q)}
    =
    \operatorname{diag}
    \left(
        \iota_{e,a}^{(q)}
    \right)_{a\in\mathcal{V}},
    \qquad
    \iota_{e,a}^{(q)}
    =
    \begin{cases}
        1, & a\in\mathcal{N}_{q}(e),\\
        0, & \text{otherwise}.
    \end{cases}
    \label{eq:local_mask}
\end{equation}
When both phase and frequency variables are included in the observation
operator, the same bus mask is applied to both components; the resulting
observation-space mask is denoted by $W_{e,H}^{(q)}$. In all results reported here, this mask is evaluated at $q=0$.

ParaEMT produces one deterministic trajectory for each faulted line and
pre-fault operating state $(e,z)$. Consequently, the covariance target used
in the calibration is not an ensemble covariance obtained from repeated EMT
realizations. Instead, it measures the scatter of the projected finite-time
EMT state across the training operating-state family at fixed faulted line.
For the global fit, $e\in\mathcal{E}_{\rm train}$ and all observations
in this target use $\Top=\Tev=0.8\,\mathrm{s}$. Let
$\mathcal{Z}_{\rm train}$ denote the 46 scaled training operating states
specified in Sec.~\ref{subsec:calibration_protocol}, and define
\begin{equation}
    \overline y^{\rm EMT}_{e}
    =
    \frac{1}{|\mathcal Z_{\rm train}|}
    \sum_{z\in\mathcal Z_{\rm train}}
    y^{\rm EMT}_{e,z}(\Tev).
    \label{eq:emt_operating_mean}
\end{equation}
The corresponding sample covariance is
\begin{equation}
\begin{aligned}
    \widehat\Sigma^{\rm EMT}_{e}
    =
    \frac{1}{|\mathcal Z_{\rm train}|-1}
    \sum_{z\in\mathcal Z_{\rm train}}&
    \left(
        y^{\rm EMT}_{e,z}(\Tev)
        -
        \overline y^{\rm EMT}_{e}
    \right)\times\\
    &\left(
        y^{\rm EMT}_{e,z}(\Tev)
        -
        \overline y^{\rm EMT}_{e}
    \right)^{\mathsf T}.
\end{aligned}    \label{eq:emt_operating_covariance}
\end{equation}
The target $\widehat\Sigma^{\rm EMT}_{e}$ measures variation across training
operating states, not numerical randomness in ParaEMT. It constrains the
variability represented by $\sigma$ in the global surrogate. The fitted
parameter distributions are reported in
Appendix~\ref{app:hierarchy_specification}.

The objective is given in \eqref{eq:calibration_objective} in the main text.
The calibration weights are fixed throughout the analysis as
\begin{equation}
    \zeta_{\rm cov}=1.0,
    \qquad
    \zeta_{\rm reg}=0.5.
    \label{eq:calibration_weights}
\end{equation}

The mean term in \eqref{eq:calibration_objective} constrains the transferred
endpoint response, while the covariance term uses the operating-state
scatter in \eqref{eq:emt_operating_covariance}. The parameters are fitted
jointly because $\sigma$ also enters the mean dynamics.

The regularization acts on the endpoint-masked predicted covariance
\begin{equation}
    \Sigma_{\Theta,e,z}^{(0)}
    =
    W_{e,H}^{(0)}
    H\Sigma_\Theta(\Tev;e,z)H^{T}
    W_{e,H}^{(0)}.
    \label{eq:predicted_local_covariance}
\end{equation}
Because the $q=0$ mask introduces zero directions outside the endpoint
subspace, the determinant is evaluated after adding the fixed diagonal
regularization
\begin{equation}
    \widetilde\Sigma_{\Theta,e,z}^{(0)}
    =
    \Sigma_{\Theta,e,z}^{(0)}
    +
    \epsilon I,
    \qquad
    \epsilon=0.01.
    \label{eq:regularized_local_covariance}
\end{equation}
The regularization term is then
\begin{equation}
    \mathcal R(\Theta)
    =
    \sum_{(e,z)\in{\cal D}_{\rm train}}
    \left|
        \log\det
        \left(
            \widetilde\Sigma_{\Theta,e,z}^{(0)}
        \right)
    \right|.
    \label{eq:logdet_regularization}
\end{equation}
The diagonal shift makes the covariance entering the determinant
nonsingular. The absolute log-determinant penalizes departures from the
reference determinant scale in either direction; it does not bound the
condition number or control each covariance eigenvalue separately. This
penalty supplements, rather than replaces, the operating-state covariance
target in \eqref{eq:emt_operating_covariance}.

The global fit uses only training-line endpoint observations at the
$0.8$-s reclosure snapshot. Models $\mathrm M_3$ and $\mathrm M_4$ additionally use
$0.8$-s observations from the assessment lines, with one pair per line
or per line--operating-state combination, respectively. Those additional
observations do not enter the globally shared fit. The coupled swing
dynamics propagate the disturbance throughout $\mathcal G$.
Network-wide comparison at $(\Top,\Tev)=(1.5,2.5)\,\mathrm{s}$ then
assesses all models at a duration and observation time not used for
identification. The faulted lines themselves are held out only for the
global models.

\subsection{Mean-Square Stability Analysis}\label{app:ms_stability}

The second-moment diagnostic characterizes the reduced stochastic dynamics.
Multiplicative noise can destabilize that moment even when the deterministic
drift is stable \cite{arnold1974stochastic}. During
$0\leq t\leq\Top$, the homogeneous operator, acting on a covariance-matrix argument $C$, is
\begin{equation}
    \mathcal{K}_{e,\Theta}(C)
    \doteq
    A_{e,\mathrm I}C
    +
    CA_{e,\mathrm I}^{\mathsf T}
    +
    \sum_{a\in\mathcal N_0(e)}
    G_{e,a}CG_{e,a}^{\mathsf T}.
    \label{eq:ms_operator}
\end{equation}
The term
$\displaystyle\sum_{a\in\mathcal N_0(e)}G_{e,a}\mu_{\Theta}\mu_{\Theta}^{\mathsf T}G_{e,a}^{\mathsf T}$
in \eqref{eq:covariance_dynamics} is a nonhomogeneous forcing term and
therefore does not enter the spectral stability criterion. Vectorization of
\eqref{eq:ms_operator} gives
\begin{equation}
\begin{aligned}
    \operatorname{vec}\!\left(
        \mathcal{K}_{e,\Theta}(C)
    \right)&
    =
    K_{e,\Theta}\operatorname{vec}(C),
    \\
    K_{e,\Theta}&
    =
    I\otimes A_{e,\mathrm I}
    +
    A_{e,\mathrm I}\otimes I
    +
    \sum_{\mathclap{a\in\mathcal N_0(e)}}
    G_{e,a}\otimes G_{e,a},
\end{aligned}
    \label{eq:ms_matrix}
\end{equation}
where $\otimes$ denotes the Kronecker product. The corresponding
mean-square stability margin is defined by
\begin{equation}
    \eta_{\rm ms}(e,z;\Theta)
    \doteq
    -
    \max_{\lambda\in\operatorname{spec}(K_{e,\Theta})}
    \operatorname{Re}(\lambda).
    \label{eq:ms_margin}
\end{equation}
Thus, $\eta_{\rm ms}>0$ implies exponential mean-square stability of the
homogeneous open-phase stochastic dynamics, whereas $\eta_{\rm ms}<0$
indicates a noise-induced second-moment instability. For a scalar
SDE $dX=aX\,dt+\sigma X\,dW$, this criterion reduces to
$2a+\sigma^2<0$, illustrating why deterministic stability alone is
insufficient in the presence of multiplicative noise
\cite{khasminskii2012stochastic}.

We evaluate \eqref{eq:ms_margin} at the globally calibrated WECC parameters
$(\alpha^*,\sigma^*)=(0.516,1.5)$ using the same inertia, damping,
network-Laplacian, and fault-local noise structure employed by the surrogate.
The calculation removes the global-phase mode before evaluating the
reported physical spectrum.
For the limiting tested contingency $X15\#228$, the rightmost physical
second-moment eigenvalue is
\[
    \max_{\lambda\in\operatorname{spec}(K_{e,\Theta})}
    \operatorname{Re}(\lambda)
    \simeq
    -0.2843~\mathrm{s}^{-1},
\]
which gives
\[
    \eta_{\rm ms}
    \simeq
    0.2843~\mathrm{s}^{-1}>0.
\]
The reported operator therefore has a positive margin at
$\sigma^{*2}=2.25~\mathrm{s}^{-1}$ for the homogeneous second moment.

Because the stochastic forcing is removed at $t=\Top$, the subsequent
interval $\Top\leq t\leq\Tev$ is governed by the deterministic covariance
operator
\begin{equation}
    \mathcal{K}_{0}(C)
    =
    AC+CA^{\mathsf T},
    \label{eq:postreclosure_ms_operator}
\end{equation}
The post-reclosure covariance therefore evolves under the nominal
deterministic operator.

\section{Supporting Calibration and Transient Diagnostics}
\label{app:validation_details}
This appendix reports the parameter distributions, the two temporal
sensitivity studies, and the full illustrative large-excursion trajectory. Their
interpretation is given with the main calibration results in
Sec.~\ref{subsec:global_wecc_calibration}.

\subsection{Parameter Variation Across Scenarios}
\label{app:hierarchy_specification}
The scenario-specific $\mathrm M_4$ fits use the additional reclosure observations
and parameter-sharing allocation described in
Sec.~\ref{subsec:calibration_protocol}. They yield 1,150 parameter pairs,
\begin{equation}
    \Theta_{e,z}
    =
    (\alpha_{e,z},\sigma_{e,z}).
\end{equation}

Across the 25 assessment lines and 46 operating states, these
$0.8$-s calibrated parameters span approximately
\begin{equation}
    0.342
    \leq
    \alpha_{e,z}^{\ast}
    \leq
    0.810,
    \qquad
    1.070
    \leq
    \sigma_{e,z}^{\ast}
    \leq
    2.643,
    \label{eq:parameter_ranges}
\end{equation}
with mean values of approximately
\begin{equation}
    \overline{\alpha}^{\ast}\approx0.558,
    \qquad
    \overline{\sigma}^{\ast}\approx1.620.
    \label{eq:parameter_means}
\end{equation}

\begin{figure}[t]
\centering
\includegraphics[width=\columnwidth]{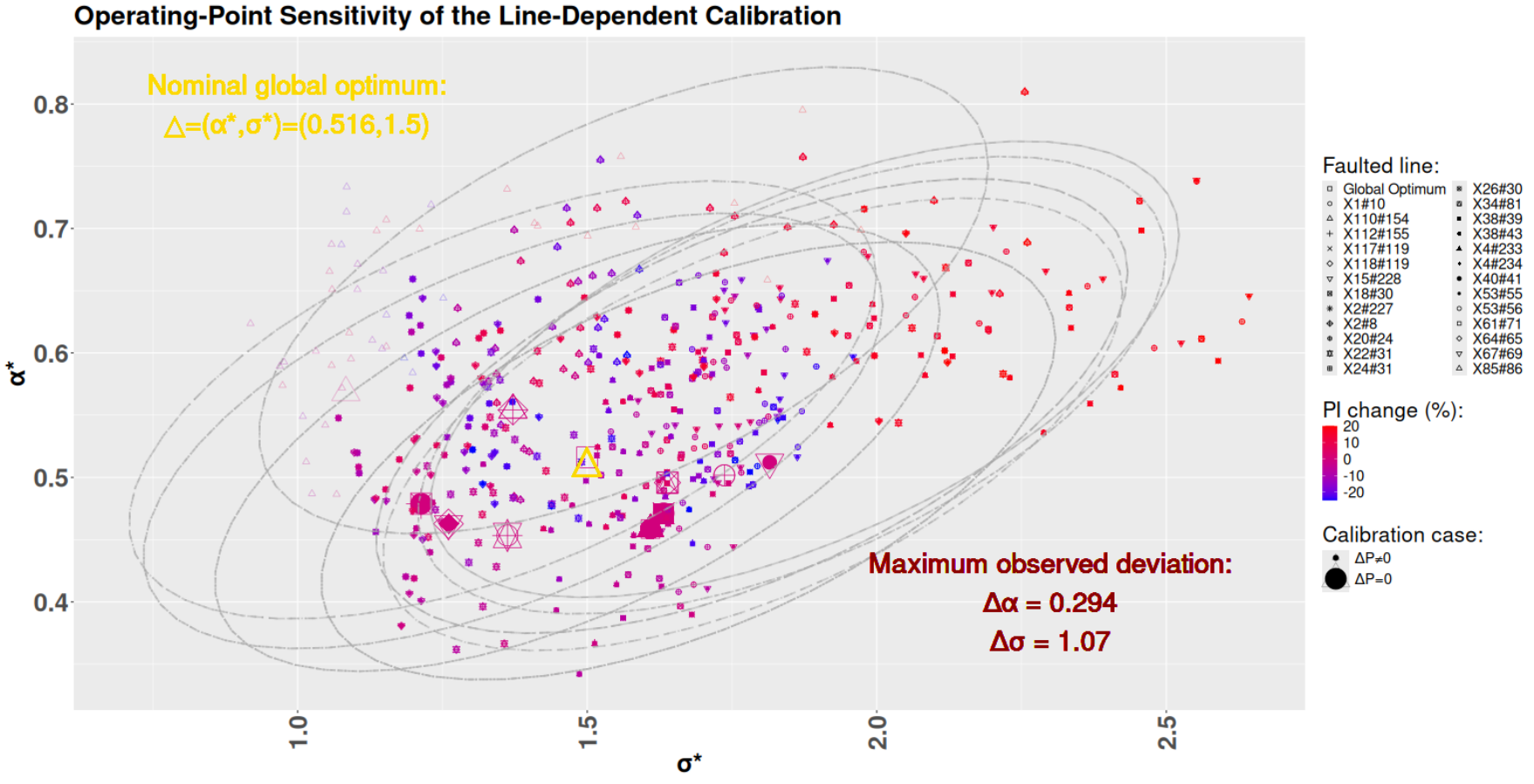}
\caption{Scenario-specific $\mathrm M_4$ parameters fitted at
$\Top=\Tev=0.8\,\mathrm{s}$ for 25 assessment lines and 46
operating states, relative to the global pair
$(\alpha^*,\sigma^*)=(0.516,1.5)$. Each of the 1,150 pairs is held
fixed in the later hierarchy assessment. Maximum reported deviations
are $\Delta\alpha=0.294$ and $\Delta\sigma=1.07$. Ellipses summarize
line-wise parameter distributions, not confidence regions.}
\label{fig:opsensitivity}
\end{figure}

Figure~\ref{fig:opsensitivity} displays the variation with fault identity
and operating condition in the respective parameter units. The
corresponding state-transfer errors and calibration-effort trade-off are
analyzed in Sec.~\ref{subsec:global_wecc_calibration}.

\subsection{Horizon Dependence of the Map Error}\label{subsec:fault_duration}

\begin{figure}[t]
\centering
\includegraphics[width=\columnwidth]{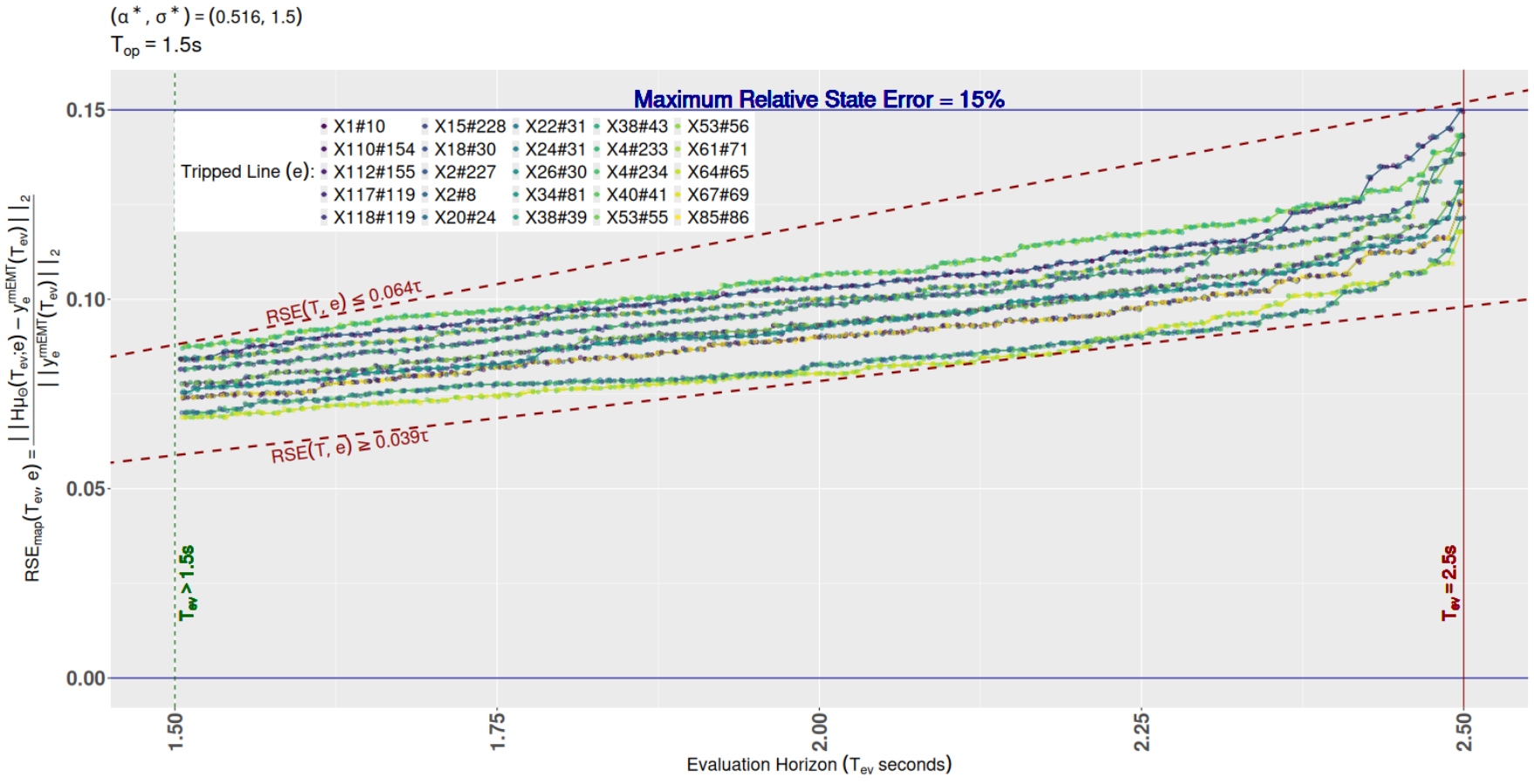}
\caption{Horizon-dependent $RSE_{\rm map}$ for the 25 WECC testing
contingencies at fixed $\Top=1.5\,\mathrm{s}$ and unscaled injections
($\chi=1$), with $(\alpha^*,\sigma^*)=(0.516,1.5)$ held fixed.
The observation time varies from reclosure to $2.5\,\mathrm{s}$.
Dashed red curves are empirical guides, not universal bounds.}
\label{fig:fault_duration}
\end{figure}

This experiment varies $\Tev$ at fixed $\Top=1.5\,\mathrm{s}$
for the 25 testing lines. The injection vector is not rescaled:
$\chi=1$ in \eqref{eq:operating_scaling}. The initial operating
condition and global pair
$\Theta^\ast=(\alpha^\ast,\sigma^\ast)=(0.516,1.5)$ remain fixed;
the pair is not refitted as the observation time changes. In the notation
below, the fixed operating state is suppressed. For the projected EMT
state $y_e^{\mathrm{EMT}}(\Tev)$, the relative map error is
\begin{equation}
    \mathrm{RSE}_{\mathrm{map}}(\Tev,e)
    =
    \frac{
        \left\|
            H\mu_{\Theta^\ast}(\Tev;e)
            -
            y_e^{\mathrm{EMT}}(\Tev)
        \right\|_2
    }{
        \left\|
            y_e^{\mathrm{EMT}}(\Tev)
        \right\|_2
    },
    \label{eq:duration_error}
\end{equation}
Figure~\ref{fig:fault_duration} shows that line-dependent discrepancies
continue to grow after reclosure, reaching approximately
$0.11$--$0.15$ at $\Tev=2.5\,\mathrm{s}$. These are nominal-injection
results, whereas Fig.~\ref{fig:boxplots} aggregates the 46 scaled
operating states. The increasing spread indicates that the subsequent
propagation retains a dependence on fault location.

\subsection{Sensitivity to Open-Phase Duration}\label{subsec:dead_time}

The complementary experiment fixes $\Tev=1.5\,\mathrm{s}$ and varies
$\Top$ from $0$ to $1.5\,\mathrm{s}$ for the same 25 testing lines.
The injections remain unscaled ($\chi=1$), and the global pair
$\Theta^\ast=(\alpha^\ast,\sigma^\ast)=(0.516,1.5)$ is not refitted.
This study changes the switching duration, not the operating-state
scale. For faulted line $e$, the error is
\begingroup
\small
\begin{equation}
\begin{aligned}
    \mathrm{RSE}_{\mathrm{map}}(\Top,e)
    &=
    \frac{
        \left\|
            H\mu_{\Theta^\ast}(\Tev;e,\Top)
            -
            y_e^{\mathrm{EMT}}(\Tev;\Top)
        \right\|_2
    }{
        \left\|
            y_e^{\mathrm{EMT}}(\Tev;\Top)
        \right\|_2
    },\\
    \Tev&=1.5\,\mathrm{s}.
\end{aligned}
    \label{eq:dead_time_error}
\end{equation}
\endgroup
Unlike \eqref{eq:duration_error}, this comparison changes the switching
time rather than the observation horizon.

\begin{figure}[t]
\centering
\includegraphics[width=\columnwidth]{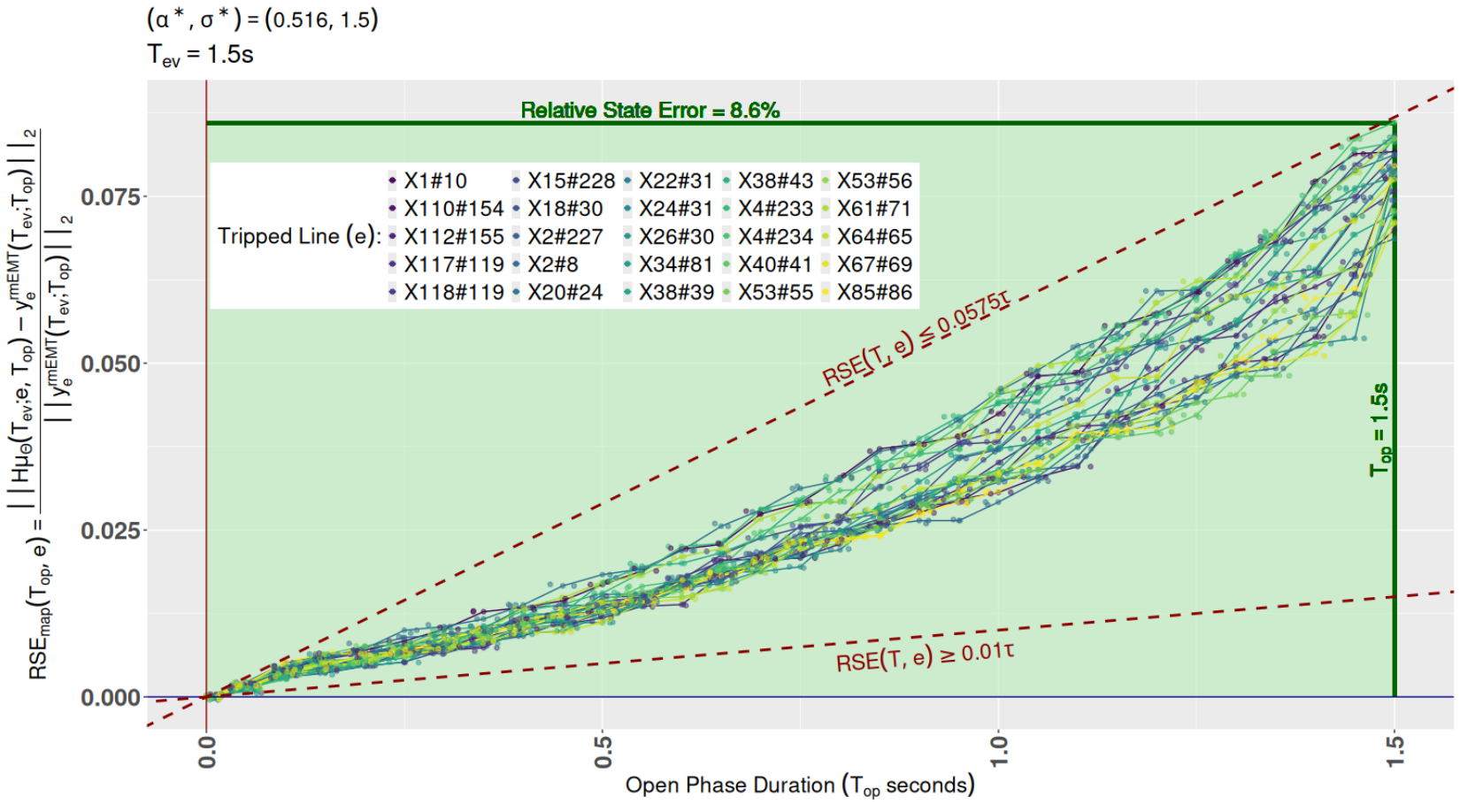}
\caption{Open-phase-duration dependence of $RSE_{\rm map}$ for the
25 WECC testing contingencies at fixed $\Tev=1.5\,\mathrm{s}$ and
unscaled injections ($\chi=1$). The global pair
$(\alpha^*,\sigma^*)=(0.516,1.5)$ is held fixed while $\Top$ varies
from $0$ to $1.5\,\mathrm{s}$. Dashed red curves are empirical guides,
not universal bounds.}
\label{fig:dead_time}
\end{figure}

Figure~\ref{fig:dead_time} shows error approaching zero as
$\Top\rightarrow0$ and increasing with the open-phase duration.
A longer interval extends the action of $A_e(\alpha^\ast)$ and the
localized multiplicative forcing in \eqref{eq:surrogate_sde}, changing
the state supplied to post-reclosure evolution. The widening spread
shows that this sensitivity depends on fault location.

The two temporal sweeps use the same random-seed configuration and agree
at $\Top=\Tev=1.5\,\mathrm{s}$, where both report
$RSE_{\rm map}\in[0.068,0.086]$. The upper endpoint is also the largest
displayed error over the tested duration range. The seed identifies the
numerical run; the exact conditional mean in
\eqref{eq:mean_dynamics} remains seed-independent for fixed inputs.

\subsection{Illustrative Transient and Large-Angle Diagnostic}\label{app:large_angle}

\begin{figure*}[t]
    \centering
    \includegraphics[width=\textwidth]{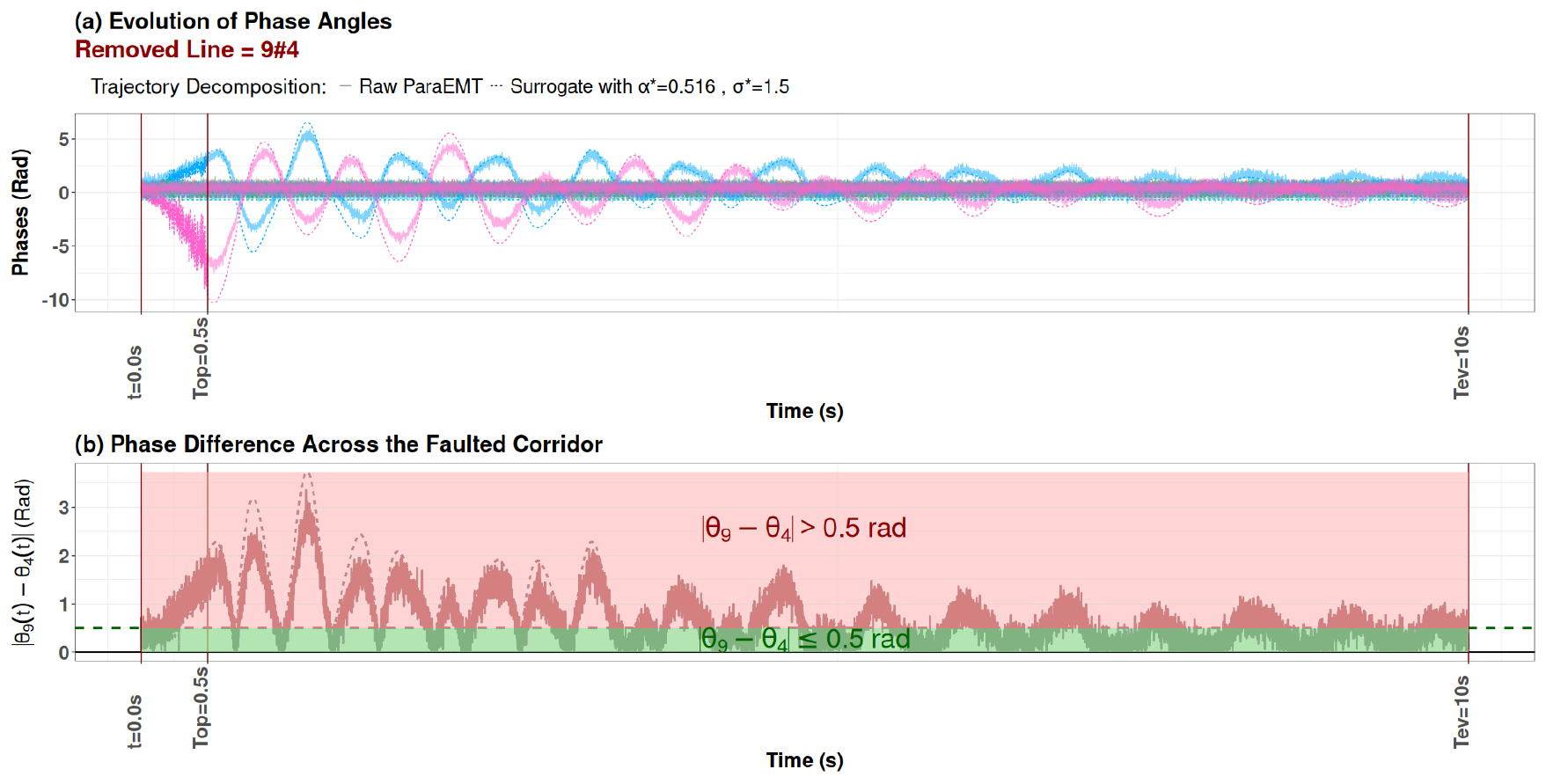}
    \caption{\textbf{Full illustrative transient and corridor-separation diagnostic.}
    The large-excursion $9\#4$ example uses $\Top=0.5\,\mathrm{s}$
    and is observed through $10\,\mathrm{s}$. The in-figure label
    $\Tev=10\,\mathrm{s}$ denotes the endpoint of this diagnostic,
    not the calibration or hierarchy-test time. Panel (a) shows phase-angle
    trajectories through reclosure; panel (b) shows
    $|\theta_9-\theta_4|$. Solid curves show ParaEMT; dashed curves show
    one surrogate realization with $(\alpha^*,\sigma^*)=(0.516,1.5)$.
    Figure~\ref{fig:open_phase_emt_surrogate} displays the open-phase
    portion of the same run. This example illustrates oscillatory
    post-reclosure behavior but does not use the $1.5$-s switching
    duration of Fig.~\ref{fig:boxplots}. Excursions beyond the shaded
    $0.5$-rad diagnostic band limit the small-angle interpretation of
    \eqref{eq:line_overload_indicator}; the band is not a protection setting.}
    \label{Point3}
\end{figure*}
The corresponding linearized no-overload region is written as
\begingroup
\small
\begin{equation}
    \Pi_{\theta}
    =
    \left\{
        \theta\in[-\pi,\pi]^{|\mathcal V|}
        :
        \left|
            \beta_{ij}
            (\theta_i-\theta_j)
        \right|
        \leq
        \bar p_{ij},
        \forall\{i,j\}\in\mathcal{E}_m
    \right\}.
    \label{eq:safe_polytope}
\end{equation}
\endgroup

Because \eqref{eq:line_overload_indicator} depends on phase differences,
uniform common-mode phase motion does not contribute directly to the
screening metric. Nevertheless, the linearized relation
$\beta_{ij}(\theta_i-\theta_j)$ assumes sufficiently small phase separation.
We examine this assumption for the illustrative large-excursion
$9\#4$ contingency in Fig.~\ref{Point3}. Its $0.5$-s open-phase
interval and $10$-s diagnostic window are independent of the calibration
and hierarchy-test timings. Panel (a) shows the absolute phase
trajectories obtained from ParaEMT and the calibrated surrogate, while
panel (b) reports the corresponding separation
$|\theta_9(t)-\theta_4(t)|$. The phase difference repeatedly exceeds the
$0.5$-rad diagnostic threshold, particularly during the large initial
transient, and remains intermittently above this level during the subsequent
response. For this illustrative case, the maximum separation is
approximately $3.77\,\mathrm{rad}$ in ParaEMT and $4.24\,\mathrm{rad}$ in
the surrogate. Thus, the large phase excursions in Fig.~\ref{Point3} cannot
be attributed solely to common-mode drift, and the small-angle approximation
is not uniformly valid throughout this severe transient.

These corridor-separation results give the diagnostic evidence for the
indicator's interpretation in Sec.~\ref{subsec:scope_limitations}.

\section{Sampling Laws, Estimators, and Statistical Records}
\label{subsec:cross_entropy}\label{app:cem_parameter_review}

We estimate the probabilities in
\eqref{eq:line_probability}--\eqref{eq:global_probability} by combining
Cross-Entropy Method (CEM) proposal adaptation with a separate
importance-sampling estimation stage
\cite{rubinstein_optimization_1997,de2005tutorial,owen2019importance}.
Both the learned line probabilities and the learned mean duration are
retained for final estimation. Adaptation favors contingencies with large
overload scores without changing the nominal probability law or refitting
the physical surrogate parameters. Throughout this appendix, $r$ denotes
the mean of the exponential proposal, in seconds, while
$\tau_0=0.1\,\mathrm{s}$ is the mean of the nominal duration law.

\runin{Sampling space and distributions.}
Let $n_{\rm E}=|\mathcal{E}|$. The contingency sampling space is
\begin{equation}
    Z=(e,\tau)
    \in
    \mathcal Z
    \doteq
    \left(\{\varnothing\}\cup\mathcal{E}\right)
    \times\mathbb R_+,
    \label{eq:cem_space}
\end{equation}
where $e$ identifies the faulted line and $\tau$ is a realization of the open-phase duration $\Top$. The screening horizon $\Thor$ is held fixed. The WECC implementation samples only faulted lines, assigning zero probability to $e=\varnothing$ under both the nominal and proposal laws. Thus, the estimates are conditional on a fault occurring within the modeled line set.

The nominal law uses uniform line selection and an exponential duration,
\begin{equation}
    p_Z(e,\tau)
    =
    \frac{1}{n_{\rm E}}\frac{1}{\tau_0}\exp\!\left(-\frac{\tau}{\tau_0}\right),
    \qquad
    \tau_0=0.1\,\mathrm{s},
    \quad e\in\mathcal{E}.
    \label{eq:cem_nominal}
\end{equation}
The adaptive proposal has categorical line probabilities and a single
shared exponential mean,
\begin{equation}
    q_{\nu}(e,\tau)
    =
    \phi_e\,\frac{1}{r} \exp\!\left(-\frac{\tau}{r}\right),
    \qquad r>0,\quad e\in\mathcal{E},
    \label{eq:cem_distribution}
\end{equation}
with $q_\nu(\varnothing,\tau)=0$. The mean $r$ is shared by all
faulted lines; no line-specific duration mean is fitted. Both $r$ and $\tau_0$ have units
$\mathrm{s}$. The nominal expected open-phase duration is therefore
$0.1\,\mathrm{s}$, whereas the proposal mean changes during adaptation.
The operating-state scaling factor in \eqref{eq:operating_scaling} is a
different quantity.

The categorical probabilities satisfy
\begin{equation}
    \phi_e\geq0,
    \qquad
    \sum_{e\in\{\varnothing\}\cup\mathcal{E}}\phi_e=1,
    \qquad \phi_{\varnothing}=0,
\end{equation}
and the adjustable parameters are
\begin{equation}
    \nu
    =
    \left((\phi_e)_{e\in\mathcal{E}},r\right).
\end{equation}
The required support condition is
\begin{equation}
    p_Z(e,\tau)>0
    \quad\Longrightarrow\quad
    q_\nu(e,\tau)>0
    \quad\text{and}\quad
    q_{\nu^\ast}(e,\tau)>0,
    \label{eq:cem_support}
\end{equation}
where $q_{\nu^\ast}$ denotes the final estimation proposal defined below.
Uniform initialization, positive duration parameters, and positive
categorical smoothing preserve support during a finite adaptation stage.

\runin{Elite selection and parameter updates.}
The WECC adaptation uses
\begin{equation}
\boxed{
\begin{gathered}
    N_{\rm CE}=250,\qquad K=20,\qquad \rho_{\rm CE}=0.1,\\
    n_{\rm E}=329,\qquad n_{\rm elite}=25,\qquad \kappa=0.9,\\
    \phi_e^{(0)}=\frac{1}{n_{\rm E}},\qquad
    r^{(0)}=1.0\,\mathrm{s}.
\end{gathered}
}
\label{eq:cem_initialization}
\end{equation}
Here, $N_{\rm CE}$ is the number of samples per adaptation iteration,
$K$ is the fixed iteration count, and $\kappa$ is the fraction of the
previous parameter value retained in each update. The symbol $\kappa$
distinguishes this numerical smoothing parameter from the physical
coupling $\alpha$.
The initial proposal mean $r^{(0)}=1.0\,\mathrm{s}$ differs
from the nominal mean $\tau_0=0.1\,\mathrm{s}$. Updating the proposal
does not change the expected duration under the nominal law.

At iteration $\ell$, draw $N_{\rm CE}$ contingencies from
$q_{\nu^{(\ell-1)}}$, evaluate the corresponding trajectories through
\eqref{eq:surrogate_sde}, and compute
$S_k^{(\ell)}=S(\theta_{0\rightarrow\Thor}^{(k)})$.
The elite set contains the indices of the 25 largest scores,
\begin{equation}
    \mathcal I_\ell
    =
    \operatorname{Top}_{25}
    \left(S_1^{(\ell)},\ldots,S_{250}^{(\ell)}\right),
    \qquad
    s_\ell^{\rm elite}=\min_{k\in\mathcal I_\ell}S_k^{(\ell)}.
    \label{eq:cem_elite_set}
\end{equation}
The operator $\operatorname{Top}_{25}$ selects exactly 25 observations;
ties do not enlarge the elite set. The empirical threshold $s_\ell^{\rm elite}$
is distinct from the event threshold $\gamma$ used for final estimation.

Elite membership is represented by the binary indicator
\begin{equation}
    b_k^{(\ell)}=\mathbb I\{k\in\mathcal I_\ell\}.
    \label{eq:cem_weights}
\end{equation}
Maximizing $\sum_k b_k^{(\ell)}\log q_\nu(e_k,\tau_k)$ gives
\begin{equation}
\begin{aligned}
    \widehat\phi_e^{(\ell)}
    &=
    \frac{1}{25}
    \sum_{k\in\mathcal I_\ell}\mathbb I\{e_k=e\},\\[1mm]
    \widehat r^{(\ell)}
    &=
    \left(
        \frac{1}{25}\sum_{k\in\mathcal I_\ell}\tau_k
    \right).
\end{aligned}
\label{eq:cem_updates}
\end{equation}
The mean-duration update uses the pooled elite durations, consistent with
the shared mean in \eqref{eq:cem_distribution}. The smoothed updates are
\begin{equation}
\begin{aligned}
    \phi_e^{(\ell)}
    &=0.9\phi_e^{(\ell-1)}+0.1\widehat\phi_e^{(\ell)},\\
    r^{(\ell)}
    &=0.9r^{(\ell-1)}+0.1\widehat r^{(\ell)}.
\end{aligned}
\label{eq:cem_smoothing}
\end{equation}
The categorical probabilities are renormalized after each update to remove
roundoff. Even when a line has no elite observations,
$\phi_e^{(\ell)}\geq0.9^{\ell}/n_{\rm E}>0$ in exact arithmetic.

This is an unweighted elite-based CE proposal-learning rule. Unlike the
nominal-density-weighted CE formulation for rare-event estimation, it does
not include $p_Z/q_\nu$ in the adaptation objective. Importance weights are
applied in the final probability estimate. Adaptation is also independent
of the chosen final threshold $\gamma$; it favors large global scores
rather than stopping when a prescribed event level is reached.

The stopping rule is a fixed computational budget,
\begin{equation}
    K=20,\qquad
    \nu^{(K)}
    =
    \left((\phi_e^{(20)})_{e\in\mathcal{E}},r^{(20)}\right).
    \label{eq:cem_stopping}
\end{equation}
No parameter-convergence tolerance or early-stopping test is used.
The largest score, elite threshold, shared mean duration, and maximum categorical
probability are recorded at every iteration.

\runin{Fresh-batch probability estimation.}
After the twentieth iteration, the complete learned proposal is frozen:
\begin{equation}
\begin{aligned}
    \phi_e^\ast&=\phi_e^{(20)},\\
    r^\ast&=r^{(20)},\\
    q_{\nu^\ast}(e,\tau)
    &=
    \phi_e^{(20)}\frac{1}{r^{(20)}}\exp\!\left(-\frac{\tau}{r^{(20)}}\right)
    =
    \phi_e^\ast \frac{1}{r^{\ast}} \exp\!\left(-\frac{\tau}{r^{\ast}}\right).
\end{aligned}
\label{eq:cem_final_proposal}
\end{equation}
Thus, $q_{\nu^\ast}=q_{\nu^{(20)}}$: the final batch retains both learned
components. The probabilities $1/n_{\rm E}$ belong to the nominal law and the
initialization, not to the final categorical proposal. In general,
$\phi_e^\ast\neq1/n_{\rm E}$. The superscript $\ast$ denotes the frozen proposal,
not a claim of global optimality or verified convergence.

Let $\xi\in\Omega_W$ denote the stochastic driving realization. Matching
the strict threshold comparison in the probability-estimation code, define
\begin{equation}
    \mathcal A_\gamma
    =
    \left\{
        (e,\tau,\xi)\in\mathcal Z\times\Omega_W:
        S\!\left(\theta_{0\rightarrow\Thor}(e,\tau,\xi)\right)
        >\gamma
    \right\},
    \label{eq:rare_event_region}
\end{equation}
for $\gamma\geq0$. The additional coordinate makes explicit that overload
is not determined by $(e,\tau)$ alone. The stochastic drivers retain the
same conditional law under the nominal and proposal simulations, so their
likelihood factors cancel. At $\gamma=0$, strict exceedance selects any
positive overload; the event $S\geq0$ would instead have probability one.

After adaptation, a fresh batch is drawn,
\begin{equation}
    (e_k,\tau_k)\overset{\mathrm{iid}}{\sim}q_{\nu^\ast},
    \qquad k=1,\ldots,N,
\end{equation}
with fresh stochastic driving realizations. Conditional on the adaptation
history, this estimation batch consists of independent draws from the frozen
proposal; no adaptation observations are reused in the estimator.

Equation~\eqref{eq:cem_estimator} in the main text gives the estimate.
The likelihood ratio corrects both changes of distribution:
\begin{equation}
\begin{aligned}
    \frac{p_Z(e_k,\tau_k)}{q_{\nu^\ast}(e_k,\tau_k)}
    &=
    \underbrace{
        \frac{1/n_{\rm E}}{\phi_{e_k}^{(20)}}
    }_{\text{line selection}}
    \underbrace{
        \frac{r^{(20)}\exp\!\left(-\frac{\tau_k}{\tau_0}\right)}
             {\tau_0\exp\!\left(-\frac{\tau_k}{r^{(20)}}\right)}
    }_{\text{open-phase duration}}\\
    &=
    \frac{r^{(20)}}
         {n_{\rm E}\phi_{e_k}^{(20)}\tau_0}
    \exp\!\left[\left(\frac{1}{r^{(20)}}-\frac{1}{\tau_0}\right)\tau_k\right].
\end{aligned}
\label{eq:cem_likelihood_ratio}
\end{equation}
The categorical factor cannot be omitted merely because the nominal fault
law is uniform. The denominator must use the probabilities that actually
generated the final line samples.

Provided that samples are generated from the density used in this ratio,
\eqref{eq:cem_support} and the fresh-batch construction imply
$\mathbb{E}[\widehat Q(\gamma)\mid\mathcal F_{\rm adapt}]=Q(\gamma)$,
where $\mathcal F_{\rm adapt}$ denotes the adaptation history.
The estimator is therefore unbiased for the specified surrogate probability
without self-normalizing the weights \cite{owen2019importance}.
This property does not remove model error relative to ParaEMT.

For weighted indicators $Y_k$ in \eqref{eq:cem_estimator}, the estimated
standard error is
\begin{equation}
    \widehat{\operatorname{se}}(\widehat Q)
    =
    \left[
        \frac{1}{N(N-1)}
        \sum_{k=1}^N(Y_k-\widehat Q)^2
    \right]^{1/2},
    \label{eq:cem_standard_error}
\end{equation}
where $Y_k$ includes the full joint likelihood ratio in
\eqref{eq:cem_likelihood_ratio}. This standard-error interpretation requires
$N>1$ and finite variance. For finitely many line categories with
$\phi_e^\ast>0$, the condition $r^\ast>\tau_0/2$ ensures finite second
moments of the weights under the mean-parameterized exponential proposal.
This follows because the duration factor in $p_Z^2/q_{\nu^\ast}$ decays only
when $2/\tau_0-1/r^\ast>0$. Small categorical probabilities can still
produce large weights; a finite empirical standard error alone does not
establish reliable tail sampling.

Algorithm~\ref{alg:n1plus} summarizes proposal adaptation and the fresh-batch
estimation stage using the complete learned joint distribution.

\begin{algorithm}[t]
\caption{N1Plus with elite-based categorical--exponential adaptation}
\label{alg:n1plus}
\begin{algorithmic}[1]

\Require Surrogate score; nominal law \eqref{eq:cem_nominal};
         screening horizon $\Thor$; threshold $\gamma\geq0$;
         final sample size $N$.

\State Set $N_{\rm CE}=250$, $K=20$, $\rho_{\rm CE}=0.1$, and $\kappa=0.9$.

\State Initialize $\phi_e^{(0)}=1/n_{\rm E}$ and
       $r^{(0)}=1.0\,\mathrm{s}$.

\For{$\ell=1,\ldots,20$}

    \State Draw $250$ contingencies from $q_{\nu^{(\ell-1)}}$
           and fresh stochastic drivers.

    \State Evaluate \eqref{eq:surrogate_sde} with $\Top=\tau_k$
           and compute the global overload scores to $\Thor$.

    \State Select the $25$ largest scores using
           \eqref{eq:cem_elite_set}.

    \State Apply
           \eqref{eq:cem_updates}--\eqref{eq:cem_smoothing}.

    \State Renormalize $\phi$ and record the adaptation diagnostics.

\EndFor

\State Freeze $r^\ast=r^{(20)}$ and
       $\phi_e^\ast=\phi_e^{(20)}$ for every $e\in\mathcal{E}$.

\State Draw $N$ fresh line indices from $\phi^{(20)}$ and durations from
       the exponential law with mean $r^{(20)}$,
       with fresh stochastic drivers.

\State Compute the global and monitored-line overload indicators.

\State Evaluate \eqref{eq:cem_estimator} using
       \eqref{eq:cem_likelihood_ratio}.

\State Report the estimate, sampling uncertainty where applicable,
       final proposal, and separate adaptation and estimation costs.

\end{algorithmic}
\end{algorithm}

The adaptation budget is $250\times20=5{,}000$ surrogate evaluations.
Consequently, a run with $N$ fresh estimation samples requires
$N+5{,}000$ evaluations; the supplied $N=35{,}000$ example therefore requires
$40{,}000$ evaluations in total. End-to-end computational comparisons must
include both stages. The supplied scorer uses $\Thor=10\,\mathrm{s}$.
When a sampled duration exceeds $\Thor$, reclosure lies beyond the screening
window; the original duration must still be retained in the likelihood ratio.

Both stages use the learned fault-location and duration components.
The categorical probabilities emphasize particular corridors, while the
shared mean duration changes the frequency of long open-phase intervals. Because
\eqref{eq:cem_distribution} factorizes, it does not introduce line--duration
dependence. Concentration of $\phi_e^{(20)}$ is an interpretable indicator of
which corridors dominate the adapted elite set, not a change in their
physical failure frequencies. Quantitative contributions to the nominal
overload probability, and conditional criticality rankings, must be
computed with the full joint likelihood ratio rather than identified with
$\phi_e^{(20)}$ itself.

The same estimator applies to a line-specific event by replacing the global
score with $S_{ij}(\theta_{0\rightarrow\Thor}^{(k)})$. A global-score proposal
can be reused when the support condition holds, although its efficiency need
not be uniform across monitored lines. Further justification and definitions
for rare-event/instanton sampling are given in Sections~VI.A--B of
\cite{almada2026realtime}; the connection to static contingency screening is
discussed in \cite{chertkov2010predicting}.

\subsection{Conditional Contingency-to-Element Estimates}\label{app:conditional_risk}

Under the uniform nominal fault law $\pi_e=1/n_{\rm E}$, a fresh CEM batch gives
\begin{equation}
    \widehat Q_{ij}(e;\gamma)
    =
    \frac{1}{N\pi_e}
    \sum_{k=1}^{N}
    w_k\,\mathbb I\{e_k=e\}
    \mathbb I\{S_{ij}^{(k)}>\gamma\},
    \label{eq:conditional_risk_estimator}
\end{equation}
where $n_{\rm E}=|\mathcal{E}|$ and $w_k$ is the full joint likelihood ratio in
\eqref{eq:cem_likelihood_ratio}. Replacing the final indicator by
$\mathbb I\{S^{(k)}>\gamma\}$ gives the global conditional estimate.
The normalization by $\pi_e$ distinguishes conditional risk from a
fault's contribution to the unconditional mixture over fault locations.
The learned proposal masses $\phi_e^{(20)}$ are not substitutes for
these weighted conditional probabilities.

\subsection{Integration of EMT-Informed Calibration}
\label{subsec:surrogate_screening_integration}

For each sampled scenario $(e,\tau)$, the dynamics in
\eqref{eq:surrogate_sde} are evaluated with $\Top=\tau$ and the global
$\mathrm M_2$ parameters
\begin{equation}
    \Theta^\ast
    =
    (\alpha^\ast,\sigma^\ast)
    =
    (0.516,1.5).
    \label{eq:screening_parameters}
\end{equation}
These physical parameters remain fixed during CE adaptation and final
estimation; they are not proposal parameters. The stochastic driving law,
pre-fault operating-state prescription, screening horizon, and overload
definitions are unchanged from the nominal assessment. Each trajectory
is scored using \eqref{eq:line_overload_indicator} and
\eqref{eq:global_overload_indicator}, and its contribution is weighted
with the complete likelihood ratio in \eqref{eq:cem_likelihood_ratio}.
The conditional normalization is given in
\eqref{eq:conditional_risk_estimator}.

\subsection{WECC Probability Diagnostics and Runtime Records}\label{app:wecc_sampling_details}

\runin{Assessment protocol.}
This appendix supplies the statistical diagnostics and timing breakdown
behind the main-text comparison in
Table~\ref{tab:reduced_data_results}. The nominal law, fixed physical
parameters, horizon $\Thor=10\,\mathrm{s}$, evaluation interval
$\Delta t=0.01\,\mathrm{s}$, and thresholds are those of
Sec.~\ref{subsec:wecc_probabilistic}. MC samples directly from the
nominal law. CEM uses the complete frozen proposal
$\phi_e^\ast=\phi_e^{(20)}$, $r^\ast=r^{(20)}$ and the likelihood
ratio in \eqref{eq:cem_likelihood_ratio}; its final batch is independent
of the adaptation draws conditional on the fitted proposal.

The reported probability estimates use the global overload event at the
four thresholds specified in the WECC study,
\begin{equation}
    Q(\gamma)
    =
    \mathbb P\!\left[
        S(\theta_{0\rightarrow\Thor})>\gamma
    \right],
    \qquad
    \gamma\in\{0,0.5,5,10\}\,\mathrm{s}.
\end{equation}

\runin{Probability estimates and sampling diagnostics.}
Table~\ref{tab:wecc_sampling_diagnostics} reports the probability estimates
and final-batch uncertainty and weight-concentration diagnostics.
Table~\ref{table4} separately records the evaluation budgets and adaptation
and estimation times. The main-text table combines only
$\gamma$, method, $N_{\rm final}$, $\widehat Q$, and total runtime;
the additional columns are retained here for reproducibility and the
variance comparison below.

For the reported estimator, define
\begin{equation}
\begin{aligned}
    w_k
    &=
    \frac{r^{(20)}}
         {n_{\rm E}\phi_{e_k}^{(20)}\tau_0}
    \exp\!\left[\left(\frac{1}{r^{(20)}}-\frac{1}{\tau_0}\right)\tau_k\right],\\
    Y_k
    &=w_k\mathbb I\{S_k>\gamma\},\\
    \widetilde w_k
    &=
    \frac{w_k}{\sum_{j=1}^{N}w_j},
    \qquad
    \widehat Q=\frac1N\sum_{k=1}^{N}Y_k,
\end{aligned}
\label{eq:wecc_weighted_samples}
\end{equation}
where $N=N_{\rm final}$,
$S_k=S(\theta_{0\rightarrow\Thor}^{(k)})$, and $n_{\rm E}=|\mathcal{E}|=329$.
The categorical factor $1/(n_{\rm E}\phi_{e_k}^{(20)})$ is retained in every
importance weight.

For independent final samples with finite-variance weighted indicators,
the estimated standard error and weight diagnostics are
\begingroup
\small
\begin{equation}
\begin{aligned}
    \widehat{\operatorname{se}}(\widehat Q)
    &=
    \left[
        \frac{\sum_{k=1}^{N}(Y_k-\widehat Q)^2}{N(N-1)}
    \right]^{1/2},\\
    \operatorname{RelSE}(\widehat Q)
    &=
    \frac{\widehat{\operatorname{se}}(\widehat Q)}{\widehat Q},\\
    \operatorname{ESS}_{w}
    &=
    \frac{\left(\sum_{k=1}^{N}w_k\right)^2}
         {\sum_{k=1}^{N}w_k^2},
    \qquad
    \widetilde w_{\max}=\max_k\widetilde w_k.
\end{aligned}
\label{eq:wecc_sampling_diagnostics}
\end{equation}
\endgroup
The relative standard error $\operatorname{RelSE}$ measures sampling
uncertainty in the probability estimate and is distinct from
$\mathrm{RSE}_{\rm map}$. The raw likelihood weight is $w_k$; its normalized diagnostic counterpart
is $\widetilde w_k$, with $\sum_k\widetilde w_k=1$.
Normalized weights are used only for diagnostics;
the probability estimate remains the ordinary importance-sampling
estimate in \eqref{eq:cem_estimator}. Weight-based ESS and maximum
normalized weight describe weight concentration, not event-specific
precision by themselves \cite{elvira2022rethinking}.

For independent, unweighted MC samples,
\begingroup
\small
\begin{equation}
    \operatorname{RelSE}_{\rm MC}
    =
    \left[
        \frac{1-\widehat Q}{(N-1)\widehat Q}
    \right]^{1/2},
    \qquad
    \operatorname{ESS}_{w}=N,
    \quad
    \widetilde w_{\max}=\frac1N.
    \label{eq:wecc_mc_relative_se}
\end{equation}
\endgroup
These relative standard errors require $\widehat Q>0$ and $N>1$.
Normalized nonnegative weights also satisfy
\begin{equation}
\begin{aligned}
    1&\leq\operatorname{ESS}_{w}\leq N,\\
    \frac1{\operatorname{ESS}_{w}}
    &\leq \widetilde w_{\max}
    \leq
    \frac{
        1+\sqrt{(N-1)(N/\operatorname{ESS}_{w}-1)}
    }{N}.
\end{aligned}
\label{eq:wecc_weight_consistency}
\end{equation}
The entries in Table~\ref{tab:wecc_sampling_diagnostics} satisfy these
identities.

\begin{table*}[t]
\caption{\textbf{Final-batch statistical reporting.}
Probability estimates, standard errors, relative standard errors,
weight-based ESS, and maximum normalized weights are obtained using the WECC simulations}
\label{tab:wecc_sampling_diagnostics}
\small
\centering
\begin{tabular}{llcccccc}
\hline
$\boldsymbol{\gamma}$ \textbf{[s]} & \textbf{Method} & $N_{\rm final}$ & $\boldsymbol{\widehat Q}$ &
$\widehat{\operatorname{se}}$ & $\operatorname{RelSE}$ &
$\operatorname{ESS}_{w}$ & $\widetilde w_{\max}$ \\
\hline
$\boldsymbol{0}$
& \textbf{CEM} & $2.500\times10^{4}$ & $\boldsymbol{5.600\times10^{-2}}$
& $1.678\times10^{-3}$ & $2.997\times10^{-2}$
& $2.009\times10^{4}$ & $5.976\times10^{-5}$ \\
& \textbf{MC} & $1.250\times10^{5}$ & $\boldsymbol{5.606\times10^{-2}}$
& $6.506\times10^{-4}$ & $1.161\times10^{-2}$
& $1.250\times10^{5}$ & $8.000\times10^{-6}$ \\
\hline
$\boldsymbol{0.5}$
& \textbf{CEM} & $4.000\times10^{4}$ & $\boldsymbol{2.628\times10^{-2}}$
& $5.430\times10^{-4}$ & $2.067\times10^{-2}$
& $3.214\times10^{4}$ & $3.817\times10^{-5}$ \\
& \textbf{MC} & $1.750\times10^{5}$ & $\boldsymbol{2.636\times10^{-2}}$
& $3.830\times10^{-4}$ & $1.453\times10^{-2}$
& $1.750\times10^{5}$ & $5.714\times10^{-6}$ \\
\hline
$\boldsymbol{5.0}$
& \textbf{CEM} & $1.050\times10^{5}$ & $\boldsymbol{7.743\times10^{-3}}$
& $2.290\times10^{-4}$ & $2.958\times10^{-2}$
& $8.437\times10^{4}$ & $1.429\times10^{-5}$ \\
& \textbf{MC} & $2.500\times10^{5}$ & $\boldsymbol{7.872\times10^{-3}}$
& $1.767\times10^{-4}$ & $2.245\times10^{-2}$
& $2.500\times10^{5}$ & $4.000\times10^{-6}$ \\
\hline
$\boldsymbol{10.0}$
& \textbf{CEM} & $1.550\times10^{5}$ & $\boldsymbol{2.723\times10^{-3}}$
& $1.940\times10^{-4}$ & $7.125\times10^{-2}$
& $1.245\times10^{5}$ & $1.573\times10^{-5}$ \\
& \textbf{MC} & $3.000\times10^{5}$ & $\boldsymbol{2.900\times10^{-3}}$
& $9.818\times10^{-5}$ & $3.385\times10^{-2}$
& $3.000\times10^{5}$ & $3.333\times10^{-6}$ \\
\hline
\end{tabular}
\end{table*}

\runin{Contribution-variance comparison.}
Comparing different sample budgets requires separating sampling precision
from the number of evaluations. Since
$N\widehat{\operatorname{se}}^{\,2}$ is the sample variance of an
individual estimator contribution, define
\begin{equation}
    \widehat V_{\rm ratio}
    =
    \frac{
        N_{\rm MC}\widehat{\operatorname{se}}_{\rm MC}^{\,2}
    }{
        N_{\rm CEM}\widehat{\operatorname{se}}_{\rm CEM}^{\,2}
    }.
    \label{eq:wecc_variance_ratio}
\end{equation}
Here, both sample sizes refer exclusively to the final estimation batches.
For independent samples, this compares the estimated per-draw variances:
values above one favor CEM, whereas values below one favor MC. It is
therefore the relevant sample-based comparison of fixed-$N$ efficiency,
rather than the ratio of the selected batch sizes. It does not account
for proposal-adaptation cost.

For the results shown in
Table~\ref{tab:wecc_sampling_diagnostics}, the ratios are
\begin{equation}
    \widehat V_{\rm ratio}
    \approx
    \begin{cases}
        0.7517, & \gamma=0,\\
        2.1761, & \gamma=0.5\,\mathrm{s},\\
        1.4184, & \gamma=5\,\mathrm{s},\\
        0.4957, & \gamma=10\,\mathrm{s}.
    \end{cases}
    \label{eq:wecc_variance_ratio_values}
\end{equation}
The contribution variance is lower for CEM at $\gamma=0.5$ and
$5\,\mathrm{s}$, and lower for MC at $\gamma=0$ and
$10\,\mathrm{s}$. The CEM final batches are smaller, so their standard
errors also reflect the different sample counts.
Equation~\eqref{eq:wecc_variance_ratio} separates that sample-size effect
from the variability of individual contributions.

\runin{Computational cost.}
The specified CE adaptation uses
\begingroup
\small
\begin{equation}
    N_{\rm adapt}=250\times20=5{,}000,
    \qquad
    N_{\rm CEM,total}=N_{\rm final}+5{,}000.
    \label{eq:wecc_sampling_budget}
\end{equation}
\endgroup
Each displayed CEM run includes one adaptation. Reusing a proposal across
thresholds incurs that adaptation cost only once. Complete CEM runtime is
\begin{equation}
    C_{\rm CEM,total}
    =
    C_{\rm adapt}+C_{\rm estimate}.
    \label{eq:wecc_wallclock_cost}
\end{equation}
This includes the pilot trajectories, sampling, elite selection,
parameter updates, final trajectory evaluations, and weighting.

Table~\ref{table4} separates adaptation time from final estimation and
records both final-batch and total evaluation counts. Adding its two time
columns gives the totals summarized in
Table~\ref{tab:reduced_data_results}; no per-realization extrapolation
is used. The comparison includes one adaptation per CEM run and pertains
to the reported sampling budgets.

\begin{table*}[t]
\caption{Separately reported CEM and MC
\cite{kroese2012monte} evaluation budgets and stage runtimes.
CEM totals include $5{,}000$ adaptation evaluations, and total time is
the sum of the supplied stage records.}
\label{table4}
\small
\centering
\begin{tabular}{llccccc}
\hline
$\boldsymbol{\gamma}$ \textbf{[s]} & \textbf{Method} & $N_{\rm final}$ & $N_{\rm total}$ &
Adaptation [s] & Estimation [s] & \textbf{Total [s]} \\
\hline
$\boldsymbol{0}$
& \textbf{CEM} & $2.500\times10^{4}$ & $3.000\times10^{4}$
& $9.241\times10^{2}$ & $3.538\times10^{3}$
& $\boldsymbol{4.462\times10^{3}}$ \\
& \textbf{MC} & $1.250\times10^{5}$ & $1.250\times10^{5}$
& $0.000$ & $2.210\times10^{4}$
& $\boldsymbol{2.210\times10^{4}}$ \\
\hline
$\boldsymbol{0.5}$
& \textbf{CEM} & $4.000\times10^{4}$ & $4.500\times10^{4}$
& $9.415\times10^{2}$ & $6.235\times10^{3}$
& $\boldsymbol{7.177\times10^{3}}$ \\
& \textbf{MC} & $1.750\times10^{5}$ & $1.750\times10^{5}$
& $0.000$ & $3.185\times10^{4}$
& $\boldsymbol{3.185\times10^{4}}$ \\
\hline
$\boldsymbol{5.0}$
& \textbf{CEM} & $1.050\times10^{5}$ & $1.100\times10^{5}$
& $9.387\times10^{2}$ & $1.876\times10^{4}$
& $\boldsymbol{1.970\times10^{4}}$ \\
& \textbf{MC} & $2.500\times10^{5}$ & $2.500\times10^{5}$
& $0.000$ & $4.536\times10^{4}$
& $\boldsymbol{4.536\times10^{4}}$ \\
\hline
$\boldsymbol{10.0}$
& \textbf{CEM} & $1.550\times10^{5}$ & $1.600\times10^{5}$
& $9.502\times10^{2}$ & $2.733\times10^{4}$
& $\boldsymbol{2.828\times10^{4}}$ \\
& \textbf{MC} & $3.000\times10^{5}$ & $3.000\times10^{5}$
& $0.000$ & $5.547\times10^{4}$
& $\boldsymbol{5.547\times10^{4}}$ \\
\hline
\end{tabular}
\end{table*}

\runin{Surrogate runtime and published ParaEMT benchmark.}
The separately reported mean computation time for one surrogate
trajectory is approximately $0.19\,\mathrm{s}$. The complete stage
timings are recorded separately in Table~\ref{table4}.

The published ParaEMT benchmark in \cite{xiong2023paraemt} reports
\begin{equation}
    c_{\rm EMT,1\ core}^{\rm pub}=29\,\mathrm{s},
    \qquad
    c_{\rm EMT,8\ cores}^{\rm pub}=28\,\mathrm{s}.
    \label{eq:wecc_emt_cost_projection}
\end{equation}
These measurements correspond to a $1$-s simulated horizon and a
$50\,\mu\mathrm{s}$ internal integration step on Windows, using a
machine equipped with two Intel Xeon Platinum 8280 2.7-GHz
processors and 512-GB RAM.

\begin{table*}[t]
\caption{Published ParaEMT runtime on the reduced WECC 240-bus
system \cite{xiong2023paraemt}. These external benchmark settings
differ from the present surrogate screening configuration.}
\label{tab:paraemt_published_benchmark}
\small
\centering
\begin{tabular}{rrrr}
\hline
Cores & Horizon [s] & EMT step [$\mu$s] & Runtime [s] \\
\hline
1 & 1 & 50 & 29 \\
8 & 1 & 50 & 28 \\
\hline
\end{tabular}
\end{table*}

The $0.01$-s trajectory-evaluation interval in the screening procedure
and the published internal EMT integration step describe different
numerical operations. Table~\ref{table4} supplies the measured stage
costs for the surrogate ensemble; Sec.~\ref{subsec:scope_limitations}
states the limits of comparisons with the external EMT benchmark.

\section{Duration-Resolved Screening}\label{app:screening_details}

\begin{figure*}[t]
\centering
\begin{subfigure}[t]{\textwidth}
\centering
\includegraphics[width=\linewidth]{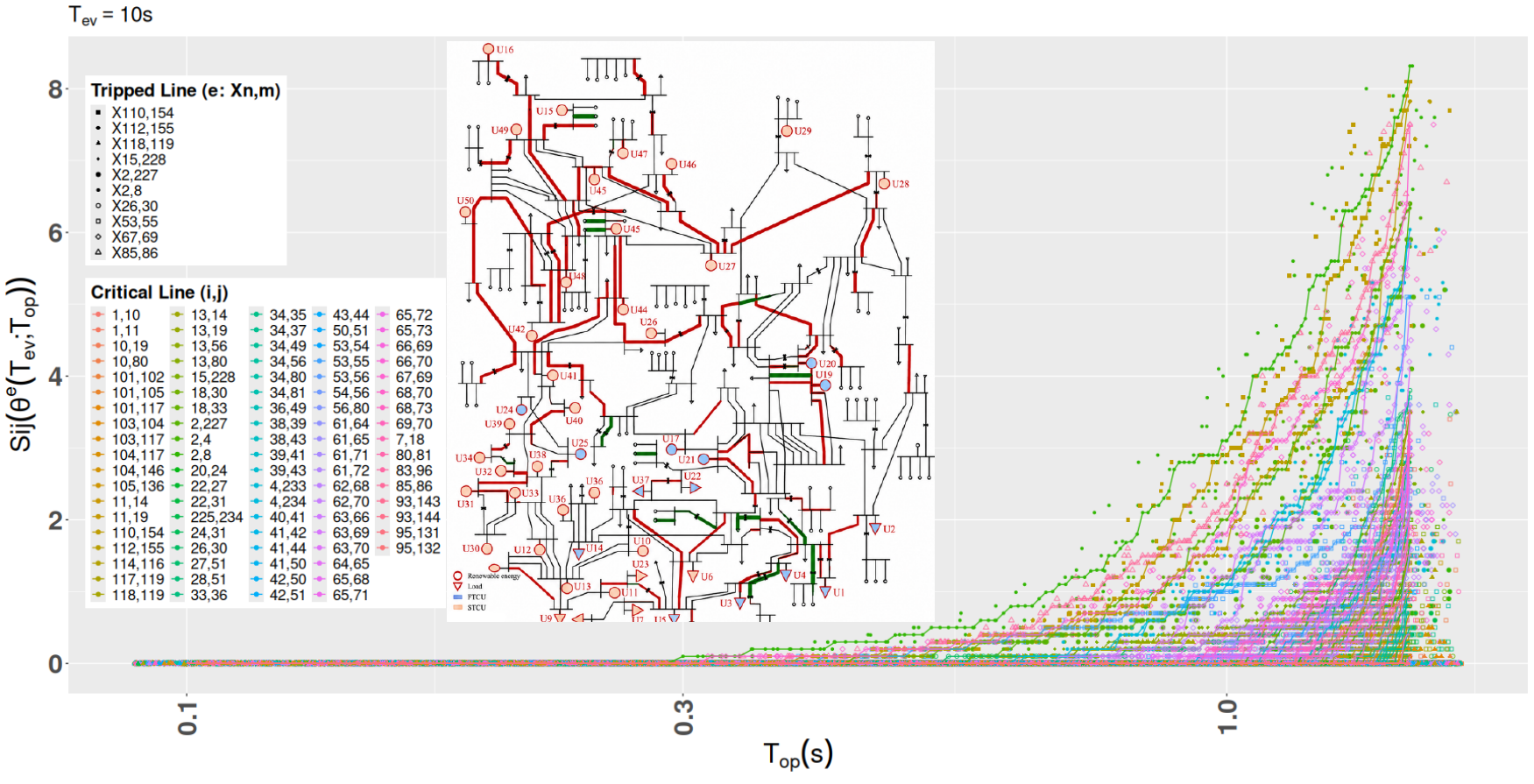}
\caption{Accumulated overload $S_{ij}$ versus $\Top$.
\textbf{\textcolor{OliveGreen}{Green}} identifies faulted lines
$(X_{n,m})$; \textbf{\textcolor{red}{red}} identifies critical
transformer-connected monitored lines $(i,j)$. The indicator measures
overload duration, not peak flow.}
\label{fig2}
\end{subfigure}
\par\medskip
\begin{subfigure}[t]{\textwidth}
\centering
\includegraphics[width=\linewidth]{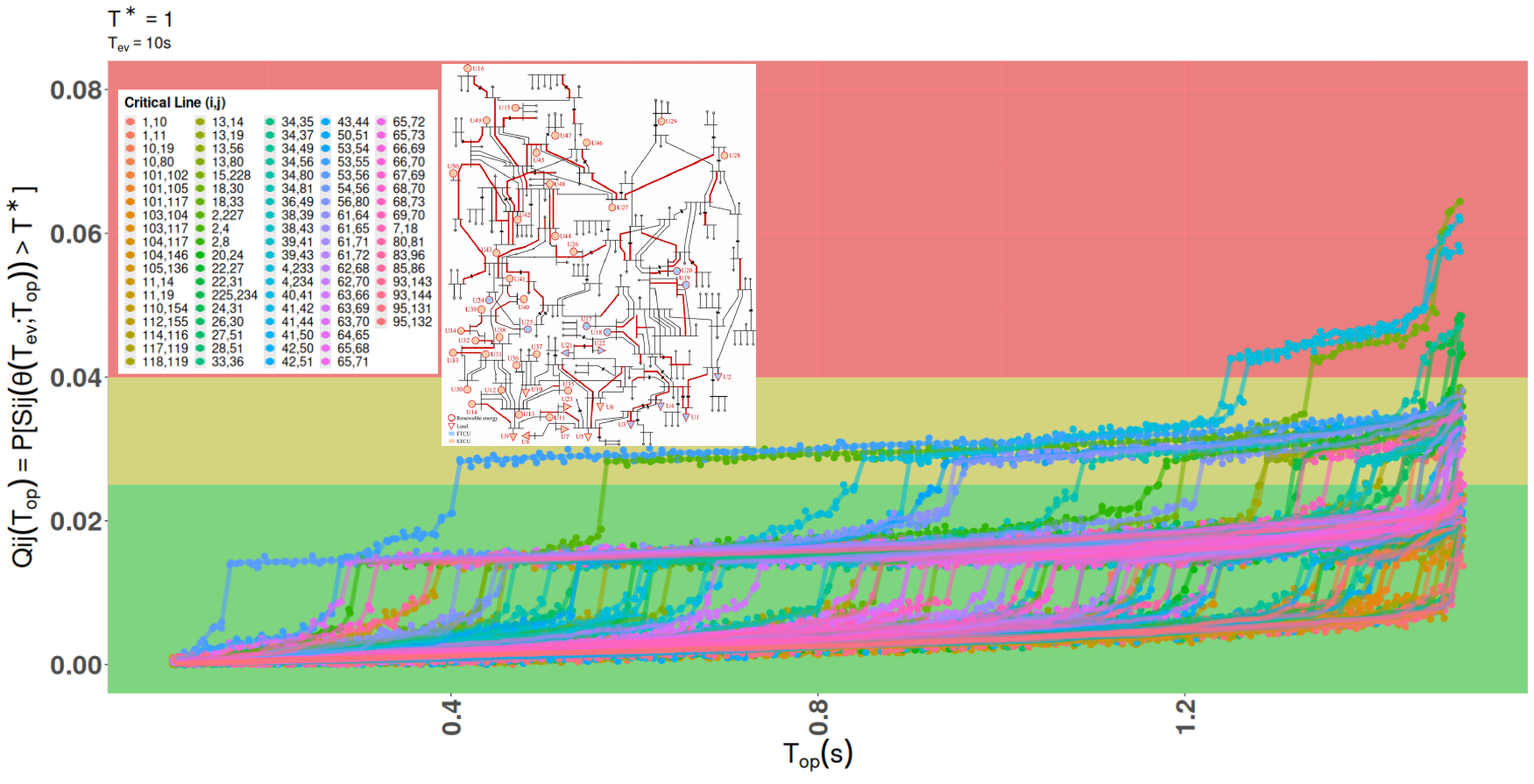}
\caption{Exceedance probability $\mathbb P[S_{ij}>\gamma]$ for
$\gamma=1\,\mathrm{s}$. Illustrative classes:
\textbf{\textcolor{OliveGreen}{safe}} ($<0.025$),
\textbf{\textcolor{YellowOrange}{warning}} ($0.025$--$0.04$), and
\textbf{\textcolor{LightRed}{emergency}} ($>0.04$).}
\label{fig3}
\end{subfigure}
\caption{Duration-resolved WECC screening: (a) accumulated line overload;
(b) line-specific exceedance probability. Both panels vary $\Top$ while accumulating exposure through
$\Thor=10\,\mathrm{s}$. The $\Tev=10\,\mathrm{s}$ notation inside
these original graphics denotes that total trajectory window, not a
calibration target. The label $T^{(*)}$ in panel~(b) denotes the
exposure threshold $\gamma=1\,\mathrm{s}$ used throughout the text. These duration-resolved quantities are distinct from
the global, duration-marginalized estimates in
Table~\ref{tab:reduced_data_results}.}
\label{fig:security_assessment}
\end{figure*}

\runin{Line-specific screening.}
Fig.~\ref{fig:security_assessment} displays line-specific screening
quantities as the open-phase duration $\Top$ varies over a fixed
$10$-s observation window. The overload integrals use the screening horizon
$\Thor$, as defined in Sec.~\ref{subsec:terminology}. The top panel
shows accumulated overload duration; the bottom panel shows the
probability of $S_{ij}>\gamma$ with $\gamma=1\,\mathrm{s}$.
This line-specific event is distinct from the global event $S>\gamma$.
The faulted element and the monitored element need not coincide.

These curves condition on duration. Their screening classes are defined
in Sec.~\ref{subsec:risk_ranking}, with the scope of their interpretation
stated in Sec.~\ref{subsec:scope_limitations}.
\runin{Relation to the main-text screening thresholds.}
The duration-resolved assessment in Fig.~\ref{fig2} identifies ten
monitored transmission lines with $S_{ij}>1\,\mathrm{s}$ for the
reported open-phase durations $\Top\geq1.0\,\mathrm{s}$; these lines
are distinct from the ten initiating faults displayed as heatmap rows.
The associated duration-resolved probabilities in Fig.~\ref{fig3}
enter the illustrative emergency class for the reported longer-duration
cases. These are duration-conditioned results for severe cases within the
tested set. The fault-conditioned, duration-marginalized probabilities in
Fig.~\ref{fig:contingency_risk_heatmap} instead lie below the first
classification boundary and use an independent continuous color scale.

\end{document}